\documentclass [11pt,a4paper]{article}

\usepackage{amsmath}
\usepackage{amsthm}
\usepackage{graphics}
\usepackage{latexsym}
\usepackage{amssymb}
\usepackage{enumerate}
\usepackage[latin2]{inputenc}
\usepackage{color}

\newtheorem  {theorem}       {Theorem}
\newtheorem  {lemma}         {Lemma}
\newtheorem  {corollary}     {Corollary}
\newtheorem  {proposition}   {Proposition}

\newtheorem* {theorem*}      {Theorem}
\newtheorem* {lemma*}        {Lemma}
\newtheorem* {corollary*}    {Corollary}
\newtheorem* {proposition*}  {Proposition}
\newtheorem* {definition*}   {Definition}
\newtheorem* {remark*}       {Remark}
\newtheorem* {remarks*}      {Remarks}

\def \N {\mathbb N}
\def \Z {\mathbb Z}
\def \Q {\mathbb Q}
\def \R {\mathbb R}

\def \T {\mathbb T}
\def \P {\mathbb P}

\def \Ordo {{\cal O}}
\def \ind{1\!\!1}

\def \lam   {\lambda}
\def \ups   {\upsilon}
\def \vphi  {\varphi}
\def \ome   {\omega}

\def \Ups   {\Upsilon}

\def\uome{\underline{\omega}}
\def\uo{\underline{\omega}}
\def\uzet{\underline{\zeta}}

\def\vareps{\varepsilon}

\def\wh{\widehat}

\def\beqs{\begin{eqnarray*}}
\def\eeqs{\end{eqnarray*}}
\def\beq{\begin{eqnarray}}
\def\eeq{\end{eqnarray}}
\def\beas{\begin{eqnarray*}}
\def\eeas{\end{eqnarray*}}
\def\bea{\begin{eqnarray}}
\def\eea{\end{eqnarray}}

\def \prob        {\ensuremath{\mathbf{P}}}
\def \expect      {\ensuremath{\mathbf{E}}}
\def \var         {\ensuremath{\mathbf{Var}}}

\def\pt{\partial_t}
\def\px{\partial_x}
\newcommand{\abs}[1]{\left|{#1}\right|}
\newcommand{\norm}[1]{\left\Vert{#1}\right\Vert}

\def \Omn {\Omega^{n}}

\def \Tn {\T^{n}}
\def \Ln {L^{n}}
\def \Kn {K^{n}}
\def \Gn {G^{n}}

\def \Xn {{\cal X}^{n}}

\def \mun {\mu^{n}}
\def \pin  {\pi^{n}}

\def \pil  {\pi^{l}}

\newcommand{\vct}[1]{ \text{\boldmath${#1}$} }

\def\vxi {\vct{\xi}}
\def\vphi{\vct{\phi}}
\def\vu  {\vct{u}}
\def\vv  {\vct{v}}
\def\vx  {\vct{x}}
\def\vlam{\vct{\lam}}

\def\vPhi{\vct{\Phi}}
\def\vvarphi{\vct{\varphi}}

\newcommand{\wih}[1] {\widehat{#1}^{n}}

\def \dom    {{\cal D}}
\def \wto    {\Rightarrow}
\def \vto    {\rightharpoonup}
\def \wvto   {\begin{array}{l}
              \displaystyle\rightharpoonup
              \\[-11.5pt]
               \hskip-0.3pt-\hskip-3.3mm\displaystyle\rightharpoonup
              \end{array}}

\title
{Derivation of the Leroux system as the hydrodynamic limit of
a two-component lattice gas}

\author {
{\sc J\'ozsef Fritz \hskip2cm Bálint Tóth}
\\[7pt]
Institute of  Mathematics
\\[3pt]
 Budapest University of Technology and Economics}

\begin{document}

\setlength{\baselineskip}{1.23\baselineskip}

\maketitle

\begin{abstract} 
The long time behavior of a couple of interacting asymmetric 
exclusion processes of opposite velocities is 
investigated in one space dimension. We do not allow two particles 
at the same site, and a collision effect (exchange) takes place
when particles of opposite velocities meet at neighboring sites.
There are two conserved quantities, and the model admits
hyperbolic (Euler) scaling; the hydrodynamic limit results in the
classical Leroux system of conservation laws, 
\emph{even beyond the appearence of shocks}. 
Actually, we prove
convergence to the set of entropy solutions, the question of
uniqueness is left open. To control rapid oscillations of Lax
entropies via 
logarithmic Sobolev inequality
estimates, the symmetric part of the
process is speeded up in a suitable way, thus a slowly vanishing
viscosity is obtained at the macroscopic level. Following
\cite{fritz1,fritz2}, 
the  stochastic version of 
Tartar--Murat theory of compensated compactness is extended to
two-component stochastic models. 
\\[7pt]
\\[7pt]
\end{abstract}
\section{Introduction}
\label{section:intro} 

The main purpose of this paper is to derive
a couple of Euler equations (hyperbolic conservation laws) in a 
regime of shocks. While the case of smooth macroscopic solutions  
is quite well understood, see 
\cite{yau} 
and 
\cite{ovy}, 
serious 
difficulties emerge when the existence of classical solutions 
breaks down. A general method to handle attractive systems has
been elaborated in 
\cite{rezakhanlou}, 
see also 
\cite{fritz1} 
and 
\cite{kipnislandim} 
for
further references. Hyperbolic models with two conservation laws,
however, can not be attractive in the usual sense because the
phase space is not ordered in a natural way. We have to extend
some advanced methods of PDE theory of hyperbolic conservation
laws to stochastic (microscopic) systems. Lax entropy and
compensated compactness are the main key words here, see
\cite{lax1}, 
\cite{lax2},  
\cite{murat}, 
\cite{tartar1},
\cite{tartar2}, 
\cite{diperna}  
for the first ideas, and the textbook 
\cite{serre} 
for a systematic treatment. The project has been initiated in 
\cite{fritz1}, 
a full 
exposition of techniques in the case of a one-component asymmetric
Ginzburg--Landau model is presented in 
\cite{fritz2}. 
Here we
investigate the simplest possible, but nontrivial two-component
lattice gas with collisions, further models are to be discussed in
a forthcoming paper 
\cite{fritztoth2}. 
Since the underlying PDE
theory is restricted to one space dimension, we also have to be
satisfied with such models. The proof is based on a strict control 
of entropy pairs at the microscopic level as prescribed by P. Lax,
L. Tartar and F. Murat for approximate solutions to hyperbolic
conservation laws. A Lax entropy is macroscopically conserved
along classical 
solutions, but the microscopic system can not have any extra
conservation law, thus we are facing with rapidly oscillating
quantities. These oscillations are to be controlled by means of
logarithmic Sobolev inequality
estimates, and effective bounds are obtainable only
if the symmetric part of the microscopic evolution is strong enough. 
That is why the 
\emph{microscopic viscosity} 
of the  model goes to infinity, i.e.  the model is changed when
we rescale it. Of course, the  
\emph{macroscopic viscosity} 
vanishes in the limit  and thus the effect
of speeding up the symmetric part of the microscopic 
infinitesimal generator  
is not seen  in the hydrodynamic limit.

Unfortunately, compensated compactness yields only
\emph{existence} 
of weak solutions, the Lax entropy condition is not sufficient for
weak 
\emph{uniqueness} 
in the case of two component systems. That is why we can 
prove convergence of the conserved fields to the set of entropy
solutions only, we do not know whether this set consists of a
single 
trajectory specified by its initial data. Let us remark that
\cite{quastelyau} 
has the same difficulty concerning the
derivation of the incompressible Navier--Stokes equation in 3
space dimensions. The Oleinik type conditions of weak uniqueness
are out of  reach of our methods because they require a one
sided uniform Lipschitz continuity of the Riemann invariants of
the macroscopic system, see 
\cite{bressan} 
for most recent results
of PDE theory in this direction. It is certainly not easy to get
such bounds at the microscopic level.

The paper is organized as follows. The microscopic model and the
macroscopic equations are introduced in the next two sections. The
main result and its conditions are formulated in Section 4. Proofs
are presented in Section 5, while some technical details are
postponed to the Appendix.

\section{Microscopic model}
\label{section:model}

\subsection{State space, conserved quantities, infinitesimal
  generator}
\label{subs:statespace}

We consider a pair of coupled asymmetric exclusion processes on
the discrete torus, particles move with an average speed 
$+1$ 
and
$-1$, 
respectively. Since we allow at most one particle per site,
the individual state space consists of three elements. There is
another effect in the interaction, something like a collision: if
two particles of opposite velocities meet at neighboring sites,
then they are also exchanged after some exponential holding times.
We can associate velocities 
$\pm1$ 
to particles according to their
categories, thus particle number and momentum are the natural
conserved quantities; the numbers of 
$+1$ 
and 
$-1$ 
particles could have been another choice.

Throughout this paper we denote by 
$\Tn$ 
the discrete torus
$\Z/n\Z$, 
$n\in\N$, 
and by 
$\T$ 
the continuous torus 
$\R/\Z$. 
The local spin space is 
$S=\{-1,0,1\}$. 
The state space of the interacting particle system of size 
$n$ 
is
\begin{equation*}
\Omn:=S^{\Tn}.
\end{equation*}
Configurations will generally be denoted as
\begin{equation*}
\uo:=(\omega_j)_{j\in \Tn}\in\Omn,
\end{equation*}
We need to separate the symmetric (reversible) part of the
dynamics. 
This  will be speeded up sufficiently in order to enhance
convergence to local equilibrium also at a mesoscopic scale.
The phenomenon of compensated compactness is materialized 
at this scale in
the hydrodynamic limiting procedure. So (somewhat artificially) we
consider separately the asymmetric and symmetric parts of the rate
functions 
$r: S\times S\rightarrow \R_+$, 
respectively, 
$s:S\times S\rightarrow \R_+$. 
The dynamics of the system consists of
elementary jumps exchanging nearest neighbor spins:
$(\omega_j,\omega_{j+1}) \rightarrow (\omega'_j,\omega'_{j+1}) =
(\omega_{j+1},\omega_{j})$, 
performed with rate 
$ \lambda r(\omega_j,\omega_{j+1}) + \kappa s(\omega_j,\omega_{j+1})$,
where
$\lambda, \kappa  > 0$ 
are  fixed speed-up  factors,  depending on
the size of the system in the limiting procedure.

The rate functions are chosen as follows:
\[
\begin{array}{rr}
 r(1,-1)=0, & \qquad r(-1,1)=2,
\\[5pt]
 r(0,-1)=0, & \qquad r(-1,0)=1,
\\[5pt]
 r(1,0)=0,  & \qquad r(0,1) =1,
\end{array}
\]
that is the rate of collisions is twice as large as that of simple
jumps, and
\[
r(\omega_j,\omega_{j+1})=\omega^-_j(1-\omega^-_{j+1})
+\omega^+_{j+1}(1-\omega^+_j)\,,
\]
where 
$\omega^+_j:=\ind_{\{\omega_j=1\}}\,,$
$\omega^-_j:=\ind_{\{\omega_j=-1\}}$ 
and 
$\ind_A$ 
denotes the
indicator of a set 
$A\,.$ 
The rates of the symmetric component are simply
\[
s(\omega_j,\omega_{j+1}) = \ind_{\{\omega_j\not=\omega_{j+1}\}}\,.
\]

The rates 
$r$ 
define a 
\emph{totally asymmetric} 
dynamics, while the rates 
$s$ 
define a 
\emph{symmetric} 
one. The infinitesimal generators defined by these rates are:
\beqs 
&& 
\Ln f(\uo) 
:= 
\sum_{j\in \Tn} 
r(\omega_j,\omega_{j+1})
(f( \Theta_{j,j+1} \uo)-f(\uo))
\\[5pt]
&& 
\Kn f(\uo) 
= 
\sum_{j\in \Tn} 
s(\omega_j,\omega_{j+1})
(f(\Theta_{j,j+1} \uo)-f(\uo))\,, 
\eeqs
where 
$\Theta_{i,j}$ 
is the spin-exchange operator,
\beqs 
\left( \Theta_{i,j} \uo \right)_{k} 
= 
\left\{
\begin{array}{lcl}
\omega_{j}    \quad   &    \text{ if }    & k=i
\\[2pt]
\omega_{i}    \quad   &    \text{ if }    & k=j
\\[2pt]
\omega_{k}    \quad   &    \text{ if }    & k\not=i,j.
\end{array}
\right. \eeqs
Recall that periodic boundary conditions are assumed in the definition 
of 
$\Ln$ 
and 
$\Kn$.

To get exactly the familiar Leroux system 
\eqref{eq:leroux} as the
limit, the two conserved quantities, 
$\eta$ 
and 
$\xi$ 
should be chosen as
\beqs
\eta_j = \eta(\omega_j) := 1-\abs{\omega_j} \quad\text{and}
\quad\xi_j = \xi(\omega_j):=\omega_j.
\eeqs
The microscopic dynamics of the model has been  defined so that
$\sum_j\xi_j$ 
and 
$\sum_j\eta_j$ 
are conserved, we shall see that
there is no room for other (independent) hidden conserved
observables. In terms of the conservative quantities we have
\begin{align}\label{eq:rates}
r(\omega_j,\omega_{j+1})&=\frac{1}{4}(1-\eta_j-\xi_j)
(1+\eta_{j+1}+\xi_{j+1})\\&+\frac{1}{4}(1+\eta_j-\xi_j)
(1-\eta_{j+1}+\xi_{j+1})\,.\notag
\end{align}

The rate functions are so chosen that the product measures
\beqs
\pin_{\rho,u}(\underline\omega)=
\prod_{j\in\Tn}\pi_{\rho,u}(\omega_j),
\eeqs
with one-dimensional marginals
\[
\pi_{\rho,u}(0)=\rho, \quad
 \pi_{\rho,u}(\pm1)=\frac{1-\rho\pm u}{2}.
\]
are stationary in time. We shall call these Gibbs measures.
The parameters take values from the set \
\[
{\cal D}:=
\{(\rho,u)\in[0,1]\times[-1,1]\,:\,\rho+\abs{u}\le1\},
\]
and the uniform 
$\pin:=\pin_{1/3,0}$ 
will serve as a reference
measure. Due to conservations, the stationary measures 
$\pin_{\rho,u}$
are not ergodic. Expectation with respect to the measures
$\pin_{\rho,u}$  
will be denoted by
$\expect_{\rho,u}(\cdot)$.
In particular, given a local observable
$\ups_i:=\ups(\omega_{i-m},\dots,\omega_{i+m})$ 
with 
$m$ 
fixed, its  equilibrium expectation will be denoted as
\beqs
\Ups(\rho,u):= \expect_{\rho,u}(\ups_i).
\eeqs

The system of microscopic size 
$n$ 
will be driven by the infinitesimal generator
\beqs
\Gn= n\Ln+ n^2\sigma\Kn,
\eeqs
where 
$\sigma=\sigma(n)$ 
is the  
\emph{macroscopic viscosity}, 
the factor 
$n\sigma(n)$ 
can be interpreted as the 
\emph{microscopic viscosity}. 
A priori we require that 
$\sigma(n)\ll1$ 
as
$n\to\infty$. 
A very important restriction,
$\sqrt{n}\sigma(n)\gg1$
will be imposed on 
$\sigma(n)$, 
see condition 
(\ref{cond:viscosity}) 
in subsection 
\ref{subs:result}.

Let 
$\mun_0$ 
be a probability distribution on 
$\Omn\,,$ 
which is the initial distribution of the microscopic system of size 
$n$,
and denote
\beqs
\mun_t:=\mun_0e^{t\Gn}
\eeqs
the distribution of
the system at (macroscopic) time 
$t$. 
The Markov process on the state space 
$\Omn$ 
driven by the infinitesimal generator 
$\Gn$,
started with initial distribution 
$\mun_0$ 
will be denoted by
$\Xn_t$.

\subsection{Fluxes}
\label{subs:fluxes}

Elementary computations show that the infinitesimal  generators
$\Ln$ 
and 
$\Kn$ 
act on the conserved quantities as follows, see
\eqref{eq:rates}.
\beqs
\begin{array}{rll}
\Ln\eta_i=
&
-
\psi(\omega_{i},\omega_{i+1})
+
\psi(\omega_{i-1},\omega_{i})
&
=:
-\psi_{i} + \psi_{i-1},
\\[5pt]
\Ln\xi_i=
&
-
\phi(\omega_{i},\omega_{i+1})
+
\phi(\omega_{i-1},\omega_{i})
&
=:
-\phi_{i} + \phi_{i-1},
\\[5pt]
\Kn\eta_i=
&
-
\psi^s(\omega_{i},\omega_{i+1})
+
\psi^s(\omega_{i-1},\omega_{i})
&
=:
-\psi^s_{i} + \psi^s_{i-1},
\\[5pt]
\Kn\xi_i=
&
-
\phi^s(\omega_{i},\omega_{i+1})
+
\phi^s(\omega_{i-1},\omega_{i})
&
=:
-\phi^s_{i} + \phi^s_{i-1},
\end{array}
\eeqs
where
\beq
\label{eq:phipsidef}
\begin{array}{rrl}
\psi_j
&=&
r(\omega_i,\omega_{i+1})\left(\eta_i-\eta_{i+1}\right)
\\[8pt]
\displaystyle
&=&
\displaystyle 
\frac12
\big\{\eta_{j}\xi_{j+1} +\eta_{j+1}\xi_{j}\big\}
+
\frac12
\big\{\eta_{j}-\eta_{j+1}\big\}
\\[8pt]
\phi_j
&=&
r(\omega_i,\omega_{i+1})\left(\xi_i-\xi_{i+1}\right)
\\[8pt]
&=& 
\displaystyle 
\frac12 \big\{ \eta_{j} + \eta_{j+1}-2 +
2\xi_{j}\xi_{j+1}\big\}+\frac12 \big\{ \xi_{j+1}\eta_{j} -
\xi_{j}\eta_{j+1} \big\} +  
\big\{\xi_{j} - \xi_{j+1}\big\}
\,,
\\[8pt]
\psi^s_j
&=&
\displaystyle
\eta_{j}-\eta_{j+1},
\\[8pt]
\phi^s_j
&=&
\displaystyle
\xi_{j}-\xi_{j+1}.
\end{array}
\eeq
Note that the microscopic fluxes of the conserved observables
induced by the symmetric rates
$s(\omega_{j},\omega_{j+1})$
are (discrete) gradients of the corresponding conserved
variables.

It is easy to compute the macroscopic fluxes:
\beq
\label{eq:macrofluxes}
\begin{array}{l}
\Psi(\rho,u)
:=
\expect_{\rho,u}(\psi_j)
=
\rho u
\\[5pt]
\Phi(\rho,u)
:=
\expect_{\rho,u}(\phi_j)
=
\rho + u^2
\end{array}
\eeq
%


\section{Leroux's equation -- a short survey}
\label{section:leroux}

Having the macroscopic fluxes (\ref{eq:macrofluxes}) computed,
the Euler equations of the system considered are expected to
be
\beq
\label{eq:leroux}
\left\{
\begin{array}{l}
\pt\rho + \px\big(\rho u\big)=0
\\[5pt]
\pt u + \px\big(\rho + u^2\big)=0.
\end{array}
\right.
\eeq
with given initial data
\beq
\label{eq:ic}
u(0,x)=u_0(x),
\quad
\rho(0,x)=\rho_0(x).
\eeq
This is exactly Leroux's equation well known in the PDE
literature, see 
\cite{serre}. 
In  the present section we shortly
review the main facts about this PDE. The first striking fact is
that such equations may have classical solutions only for some
special initial data, in general shocks are developed in a finite
time. Therefore solutions should be understood in a weak
(distributional) sense, and there are many weak solutions for the
same initial values.

The following vectorial notations sometimes make our formulas more
compact:
\beqs
\vu:=
\left(
\begin{array}{c}
\rho 
\\[5pt] 
u
\end{array}
\right),
\quad
\vPhi:=
\left(
\begin{array}{c}
\Psi 
\\[5pt] 
\Phi
\end{array}
\right),
\eeqs
\beqs
\nabla:=
\left(
\begin{array}{cc}
\displaystyle
\frac{\partial}{\partial\rho}&  
\displaystyle
\frac{\partial}{\partial u}
\end{array}
\right),
\quad
\quad
\nabla^2:=
\left(
\begin{array}{cc}
\displaystyle
\frac{\partial^2}{\partial\rho^2}&
\displaystyle
\frac{\partial^2}{\partial\rho\partial u}
\\[8pt]
\displaystyle
\frac{\partial^2}{\partial\rho\partial u}&
\displaystyle
\frac{\partial^2}{\partial u^2}
\end{array}
\right)
\eeqs
We shall use  alternatively, at convenience, the compact
vectorial and the explicit notation.

\subsection{Lax entropy pairs}
\label{subs:lent}

In the case of classical solutions 
(\ref{eq:leroux}) 
can be written as 
$\partial_t\vu+D(\vu)\partial_x\vu=0,$ 
where
\beqs
D(\rho,u):=\nabla\vPhi(\rho,u) = \left(
\begin{array}{cc}
u  & \rho
\\[5pt]
1  & 2 u
\end{array}
\right) \eeqs is the matrix of the linearized system.  The
eigenvalues of 
$D$ 
are just
\beqs 
&& 
\lambda 
= 
\lambda(\rho,u) 
:=
u+ \frac12 \left\{ \sqrt{u^2 + 4\rho} + u \right\}\,,
\\
&& 
\mu 
= 
\mu(\rho,u) 
:= 
u- \frac12 \left\{ \sqrt{u^2 + 4\rho} -u\right\}\,.
\eeqs
This means that 
(\ref{eq:leroux}) 
is
\emph{strictly hyperbolic} 
in the domain
\beqs 
\{(\rho,u):
\rho\ge0, u\in\R, (\rho,u)\not=(0,0)\}\,,
\eeqs
with marginal degeneracy (i.e. coincidence of the two characteristic
speeds, 
$\lambda=\mu$) 
at the point 
$(\rho,u)=(0,0)$.

\emph{Lax entropy/flux pairs} 
$\bigl(S(\vu),F(\vu)\bigr)$ 
are solutions of the
linear hyperbolic system 
$\nabla F(\vu)=\nabla S(\vu)\cdot\nabla\vPhi(\vu)\,,$ 
that is 
$\partial_tS(\vu)+\partial_xF(\vu)=0$ 
along classical solutions. This means that an entropy 
$S$ 
is a conserved observable. In our particular case this reads
\beq
\label{entropia}
\left\{
\begin{array}{ll}
F'_\rho=&uS'_\rho +S'_u,
\\[5pt]
F'_u=&\rho S'_\rho +2uS'_u.
\end{array}
\right.
\eeq
or, written as a second order linear equation for 
$S$:
\beq
\label{entropia2}
\rho S''_{\rho\rho}+u S''_{\rho u} -S''_{uu}=0.
\eeq
This equation is known to have many convex solutions, see
\cite{lax1}. 
We call an entropy/flux pair 
\emph{convex} 
if the map
$(\rho,u)\mapsto S(\rho,u)$ 
is convex. In particular, a globally
convex Lax entropy/flux pair defined on the whole half plane 
$\R_+\times\R$
is
\[
S(\rho,u):=\rho\log\rho +\frac{u^2}{2},
\quad
F(\rho,u):=u\rho+u\rho\log\rho+\frac{2u^3}{3}\,.
\]

Weak solutions of (\ref{entropia}) 
are called \emph{generalized entropy/flux pairs}. 
Riemann's method of solving second order linear hyperbolic PDEs
in two variables (see Chapter 4 of \cite{john}) and compactness
of $\dom$ imply that generalized entropy/flux pairs can be
approximated pointwise 
by twice differentiable entropy/flux
pairs. 

An 
\emph{entropy solution} 
of the Cauchy problem
(\ref{eq:leroux}), 
(\ref{eq:ic}) 
is a measurable function
$[0,T]\times\T\ni(t,x)\mapsto \vu(t,x)\in\R_+\times\R$ 
which for any convex entropy/flux pair 
$(S,F),$ 
and any nonnegative test function 
$\varphi:[0,T]\times\T\to\R$ 
with support in
$[0,T)\times\T$ 
satisfies
\beq
\notag
&&
\int_0^t\int_\T
\left(
\pt\varphi(t,x) S(\vu(t,x)) +
\px\varphi(t,x) F(\vu(t,x))
\right)
\,dx\,dt
\\
\label{eq:wlax}
&&
\hskip2cm
+
\int_\T
\varphi(0,x)S(\vu(0,x))
\,dx
\ge0
\eeq
Note that 
$S(\rho,u)=\pm\rho, F(\rho,u)=\pm\rho u$, 
respectively,
$S(\rho,u)=\pm u, F(\rho,u)=\pm(\rho + u^2)$ 
are entropy/flux
pairs, thus entropy solutions are (a special class of) weak
solutions. Entropy solutions of the Cauchy problem
(\ref{eq:leroux}), 
(\ref{eq:ic}) 
form a (strongly) closed subset of the Lebesgue space 
$L^p([0,T]\times\T, \,dt\,dx)=:L^p_{t,x}$
for any 
$p\in[1,\infty)$.

\subsection{Young measures, measure valued entropy solutions}
\label{subs:young}

A Young measure on
$([0,T]\times\T)\times{\cal D}$
is
$\nu=\nu(t,x;d\vv)$,
where
\\
(1)
for any
$(t,x)\in [0,T]\times\T$
fixed,
$\nu(t,x;d\vv)$
is a probability measure on
${\cal D}$,
and,
\\
(2)
for any
$A\subset{\cal D}$
fixed the map
$(t,x)\mapsto\nu(t,x;A)$
is measurable.
\\
Given a probability measure $\nu$ on $\R_+\times\R$, we shall use the
notation  
\[
\langle\nu\,,\,f\rangle:=
\int_{\cal D} f(\vv)\,\nu(d\vv).
\]
The set of Young measures will be denoted by
${\cal Y}$.
A sequence
$\nu^{n}\in{\cal Y}$
\emph{converges vaguely} to
$\nu\in{\cal Y}$,
denoted
$\nu^n\vto\nu$,
if  for any
$f\in C([0,T]\times\T\times{\cal D})$
\beqs
\lim_{n\to\infty}
\int_0^T\int_\T
\langle\nu^n(t,x)\,,\,f(t,x,\cdot)\rangle\,dt\,dx
=
\int_0^T\int_\T
\langle\nu(t,x)\,,\,f(t,x,\cdot)\rangle\,dt\,dx,
\eeqs
or, equivalently, if for any test function
$\varphi\in C([0,T]\times\T)$
and any
$g\in C({\cal D})$
\beqs
\lim_{n\to\infty}
\int_0^T\int_\T
\varphi(t,x)
\langle\nu^n(t,x)\,,\,g\rangle
\,dt\,dx
=
\int_0^T\int_\T
\varphi(t,x)
\langle\nu(t,x)\,,\,g\rangle
\,dt\,dx.
\eeqs
The set
${\cal Y}$
of Young measures will be endowed with the vague
topology induced by this notion of convergence.
${\cal Y}$
endowed
with the vague topology is  metrizable, separable and compact.
We also consider (without explicitly denoting this) the Borel
structure on 
$\cal Y$, 
induced by the vague topology.

We say that the Young measure
$\nu(t,x;d\vv)$
is 
\emph{Dirac-type} 
if there exists a measurable function
$\vu:[0,T]\times\T\to{\cal D}$
such that for almost all
$(t,x)\in[0,T]\times\T$, 
$\nu(t,x;d\vv)=\delta_{\vu(t,x)}(d\vv)$.
We denote the subset of Dirac-type  Young measures by
${\cal U}\subset{\cal Y}$.
It is a fact (see Chapter 9 of 
\cite{serre}) 
that
\[
{\cal Y}
=
\overline{\text{co}({\cal Y})}
=
\overline{\text{co}({\cal U})}
=
\overline{{\cal U}},
\]
where `co' stands for convex hull and closure is meant according to
the vague topology.

We say that the Young measure
$\nu(t,x;d\vv)$
is a 
\emph{measure valued entropy  solution} 
of the  Cauchy problem  
(\ref{eq:leroux}), 
(\ref{eq:ic})  
iff for any convex entropy/flux pair
$(S,F)$
and any positive test function
$\varphi:[0,T]\times\T\to\R_+$
with support in
$[0,T)\times\T$,
\beq
\notag
\int_0^T\int_\T
\big(
\pt\varphi(t,x)
\langle\nu(t,x)\,,\,S\rangle 
+ 
\px\varphi(t,x)
\langle\nu(t,x)\,,\,F\rangle
\big)
\,dx\,dt
\qquad\qquad
&&
\\[5pt]
\label{eq:mlax}
\hskip2cm
+
\int_\T
\varphi(0,x)
\langle\nu(0,x)\,,\,S\rangle
\, dx
\ge
0
&&
\eeq
holds true. Measure valued entropy solutions of the Cauchy problem 
(\ref{eq:leroux}), 
(\ref{eq:ic}) 
form a (vaguely) 
\emph{closed}
subset of
$\cal Y$.

Clearly, if
$\vu:[0,T]\times\T\to{\cal D}$
is an entropy  solution of the Cauchy problem
(\ref{eq:leroux}), 
(\ref{eq:ic}) 
in the sense of 
(\ref{eq:wlax}),
then the Dirac-type Young measure
$\nu(t,x;d\vv):=\delta_{\vu(t,x)}(d\vv)$
is a measure valued entropy solution in the sense of
(\ref{eq:mlax}).
The convergence of subsequences of approximate solutions to measure
solutions is  almost immediate by vague compactness, the crucial issue
is to show the Dirac property of measure valued entropy
solutions. This is the aim of the theory  of compensated compactness.

\subsection{Tartar factorization}
\label{subs:tarfac}

A probability measure $\nu(d\rho,du)$ on $\R^2$ satisfies
the \emph{Tartar factorization} property with respect to a couple
$(S_i,F_i)\,,$ $i=1,2$ of entropy/flux pairs if
\begin{equation}\label{tarfac}
\langle\nu,S_1F_2-S_2F_1\rangle=\langle\nu,S_1\rangle
\langle\nu,F_2\rangle-\langle\nu,S_2\rangle\langle\nu,F_1\rangle\,.
\end{equation}
Dirac measures certainly posses this property, and in some cases,
there is a converse statement, too.
 
The following one-parameter families of entropy/flux pairs play an
essential role in the forthcoming argument:
\beq
\label{parent}
\begin{array}{ll}
S_a(\rho,u):=\rho+au-a^2\,,
&
F_a(\rho,u):=(a+u)S_a(\rho,u)\,,
\\[5pt]
\Bar{S}_a(\rho,u):=|\rho+au-a^2|\,,
&
\Bar{F}_a(\rho,u):=(a+u)\Bar{S}_a(\rho,u)\,,
\end{array}
\eeq
where the parameter, $a\in\R\,.$ The case of $(S_a,F_a)$ is obvious
because it is  a linear function of the basic conserved observables and
their fluxes.

The pair $(\Bar{S}_a,\Bar{F}_a)$ satisfies
\eqref{entropia} in the generalized (weak) sense. 
This is due to the facts that 
the line of non-differentiability, 
$\rho+au-a^2=0$, is just a characteristic line of the
PDE (\ref{entropia}), and  
$(\Bar{S}_a,\Bar{F}_a)$ coincides with $(\pm S_a,\pm F_a)$ on the domains
$D_{\pm}:=\{\pm(\rho+au-a^2)>0\}$.
 
\begin{lemma}
\label{lemma:dirac}
Suppose that a compactly supported probability measure, $\nu$ on
$\R^2$ satisfies \eqref{tarfac} for any two entropy/flux pairs
of type \eqref{parent}. Then $\nu$ is concentrated to a single
point, i.e. it is a Dirac mass.
\end{lemma}

 
\begin{proof} This is Exercise 9.1 in \cite{serre}, where detailed
instructions are also added. For Reader's convenience we reproduce
the easy proof. Suppose first that $S_a=\rho+au-a^2=0$ $\nu$-a.s.
for some $a\in\R\,,$ then
$\langle\nu,\rho\rangle+a\langle\nu,u\rangle=a^2\,;$ let $a_1$ and
$a_2$ denote the roots of this equation. Since $S_{a_1}(\rho,u)=0$
implies $S_{a_2}(\rho,u)=0,$ $u=a_1+a_2$ and $\rho=-a_1a_2$
$\nu$-a.s. Therefore we may, and do assume that
\[
g(a):=\frac{\langle\nu,\Bar{F}_a\rangle}
{\langle\nu,\Bar{S}_a\rangle}
=a+\frac{\langle\nu,u|\rho+au-a^2|\rangle}
{\langle\nu,|\rho+au-a^2|\rangle}
\]
is well defined for all $a\in\R\,.$ It is plain that $g(a)-a$ is
continuous, bounded, and $g(a)-a\to\langle\nu,u\rangle$ as
$a\to\pm\infty.$
 
Applying \eqref{tarfac} to $(S_a,F_a)$ and $(\Bar{S}_a,\Bar{F}_a)$
we get $g(a)-a=\langle\nu,uS_a\rangle/\langle\nu,S_a\rangle\,.$ On
the other hand, from \eqref{tarfac} for $(S_a,F_a)$ and
$(S_b,F_b)$ we get
\begin{equation}\label{ba}
(b-a)\left(\langle\nu,S_aS_b\rangle
-\langle\nu,S_a\rangle\langle\nu,S_b\rangle\right)
=\langle\nu,S_a\rangle\langle\nu,uS_b\rangle
-\langle\nu,S_b\rangle\langle\nu,uS_a\rangle\,,
\end{equation}
 thus dividing by
 $(b-a)\langle\nu,S_a\rangle\langle\nu,S_b\rangle\,,$ and letting
$b\to a$ we see that $g$ is differentiable, and $g'(a)\geq1,$
consequently $g(a)=a+\langle\nu,u\rangle$ for all $a\in\R\,.$ This
means that
\[
\langle\nu,\rho u\rangle+a\langle\nu,u^2\rangle
-a^2\langle\nu,u\rangle=\langle\nu,\rho\rangle
\langle\nu,u\rangle+a\langle\nu,u\rangle^2 -a^2\langle\nu,u\rangle
\]
for all $a\in\R\,,$ whence $\nu(u^2)=\nu^2(u)\,.$ Substitute now
$u=\langle\nu,u\rangle$ back into \eqref{ba}. Since
$b-a\neq\langle\nu,u\rangle$ may be assumed, we have
\[
\langle\nu,S_aS_b\rangle
=\langle\nu,S_a\rangle\langle\nu,S_b\rangle\,,
\]
consequently
$\langle\nu,\rho^2\rangle=\langle\nu,\rho\rangle^2\,.$
\end{proof}
 
This lemma establishes that measure solutions satisfying Tartar's
factorization property \eqref{tarfac} are, in fact, weak
solutions.
 

\section{The hydrodynamic limit under Eulerian scaling}
\label{section:hdl}

\subsection{Block averages}
\label{subs:blocks}

We choose a 
\emph{mesoscopic} 
block size
$l=l(n)$.
A priori
\beqs
1
\ll
l(n)
\ll
n,
\eeqs
but more serious restrictions will be imposed, see condition
(\ref{cond:blocksize}) 
in subsection 
\ref{subs:result}.
and define the 
\emph{block averages} 
of local observables in the following way:
\\
We fix once for ever a weight function
$a:\R\to\R_+$.
It is assumed that:
\\
(1)
$x\mapsto a(x)$
has support in the compact interval
$[-1,1]$,
\\
(2)
it has total weight
$\int a(x)\,dx=1$,
\\
(3)
it is even:
$a(-x)=a(x)$,
and
\\
(4)
it is twice continuously differentiable.

Given a local variable
$\ups_i$
its block average 
\emph{at macroscopic space} 
$x$ 
is defined as
\beq
\label{eq:blocks}
\wih{\ups}(x)
=
\wih{\ups}(\uo,x)
:=
\frac{1}{l}\sum_ja\left(\frac{nx-j}{l}\right)\ups _j.
\eeq
Note that, since
$l=l(n)$,
we do not denote explicitly dependence of
the block average on the mesoscopic block size
$l$.

We shall use the handy (but slightly abused) notation
\beqs
\wih{\ups}(t,x):=\wih{\ups}(\Xn_t,x).
\eeqs
This is the empirical block average process of the local observable
$\ups_i$.

In accordance with the compact vectorial notation introduced at the
beginning of Section 
\ref{section:leroux} 
we shall denote
\beqs
\vxi_j
:=
\left(\!\!
\begin{array}{c}
\eta_j
\\
\xi_j
\end{array}
\!\!\right),
\,\,\,
\vphi_j
:=
\left(\!\!
\begin{array}{c}
\psi_j
\\
\phi_j
\end{array}
\!\!\right),
\,\,\,
\wih{\vxi}(x)
:=
\left(\!\!
\begin{array}{c}
\wih{\eta}(x)
\\
\wih{\xi}(x)
\end{array}
\!\!\right),
\,\,\,
\wih{\vphi}(x)
:=
\left(\!\!
\begin{array}{c}
\wih{\psi}(x)
\\
\wih{\phi}(x)
\end{array}
\!\!\right),
\eeqs
and so on.

Let
$\wih{\vxi}(t,x)$
be the sequence of empirical block average
processes of the
conserved quantities, as defined above, regarded as elements of
$L^1_{t,x}:=L^1([0,T]\times\T)$.
We denote by 
$\P^n$ 
the distribution of these in
$L^1_{t,x}$:
\beq
\label{eq:l1distr}
\P^n(A)
:=
\prob\left(\vxi^n\in A\right),
\eeq
where
$A\in L^1_{t,x}$
is (strongly) measurable. Tightness and weak
convergence of the sequence of probability measures
$\P^n$
will be meant  according to  the norm (strong) topology of
$L^1_{t,x}$.
Weak convergence of a subsequence
$\P^{n'}$
will be denoted
$\P^{n'}\wto\P$.

Further on, we denote by
$\nu^n$
the sequence of Dirac-type random Young
measures  concentrated on the trajectories of the empirical averages
$\wih{\vxi}(t,x)$
and by
$\Q^n$
their distributions on
${\cal Y}$:
\beq
\label{eq:young}
\nu^n(t,x;d\vv)
:=
\delta_{\wih{\vxi}(t,x)}(d\vv),
\qquad
\Q^n(A)
:=
\prob\left(\nu^n\in A\right),
\eeq
where
$A\in{\cal Y}$
is (vaguely)  measurable.
Due to vague compactness of
$\cal Y$,
the sequence of probability
measures
$\Q^n$
is automatically tight.  Weak convergence of a subsequence
$\Q^{n'}$
will be meant according  to the vague topology of
$\cal Y$
and will be denoted
$\Q^n\wvto\Q$.
In this case we shall also say that the subsequence of random Young
measures  
$\nu^{n'}$ 
(distributed according to
$\Q^{n'}$)
\emph{converges vaguely in distribution} 
to the random Young measure
$\nu$
(distributed according to
$\Q$),
also denoted
$\nu^n\wvto\nu$.

\subsection{Main result}
\label{subs:result}

All results are valid under the following conditions

\begin{enumerate}[(A)]

\item
\label{cond:viscosity}
The macroscopic viscosity
$\sigma=\sigma(n)$
satisfies
\beqs
\label{eq:viscosity}
n^{-1/2}\ll\sigma\ll1.
\eeqs

\item
\label{cond:blocksize}
The mesoscopic block size
$l=l(n)$
is chosen so that
\beqs
\label{eq:blocksize}
n^{2/3}\sigma^{1/3}
\ll
l
\ll
n\sigma
\eeqs

\item
\label{cond:iniprofile}
The initial density profiles converge weakly in probability (or,
equivalently in any
$L^p$, $1\le p<\infty$).
That is:
for any test function
$\vvarphi:\T\to\R\times\R$
\beqs
\label{eq:iniprofile}
\lim_{n\to\infty}\expect
\Big(
\abs{
\int_\T \vvarphi(x)\cdot\big(\wih{\vxi}(0,x)-\vu_0(x)\big)
\,dx
}
\Big)
=0.
\eeqs
\end{enumerate}

Our main result is the following

\begin{theorem}
\label{thm:main}
Conditions 
(\ref{cond:viscosity}), 
(\ref{cond:blocksize}),
and 
(\ref{cond:iniprofile}) 
are in force. The sequence of probability measures
$\P^n$
on
$L^1_{t,x}$,
defined in 
(\ref{eq:l1distr}) 
is tight (according to the norm topology of
$L^1_{t,x}$).
Moreover, if
$\P^{n'}$
is a   subsequence
which converges weakly (according to the norm topology of
$L^1_{t,x}$), 
$\P^{n'}\wto\P$,
then the limit probability measure
$\P$
is concentrated on the entropy
solutions of the Cauchy problem  
(\ref{eq:leroux}), 
(\ref{eq:ic}).
\end{theorem}

\noindent
{\bf Remark:}
Assuming 
\emph{uniqueness} 
of the entropy solution
$\vu(t,x)$
of the Cauchy  problem 
(\ref{eq:leroux}), 
(\ref{eq:ic}), 
we could conclude that
\[
\wih{\vxi}\buildrel{L^1_{t,x}}\over{\longrightarrow}\vu,
\qquad
\text{ in probability.}
\]

\section{Proof}
\label{section:proof}

\subsection{Outline of proof}
\label{subs:notations}

We broke up the proof into several subsections according to what
we think to be a logical and transparent structure.

In subsection
\ref{subs:apriori} 
we state  the precise 
\emph{quantitative form} 
of the convergence to local equilibrium: 
the logarithmic Sobolev ineqaulity valid for our model and Varadhan's
large deviation bound on space-time averages of block variables. 
As main consequence of these we obtain our a priori estimates: 
the so-called 
\emph{one-block estimate}
and a  version of the so-called 
\emph{two-block estimate}, 
formulated for spatial derivatives of the empirical block
averages. These  estimates are of course the main 
probabilistic ingredients of the
further  arguments. The proof of these estimates is postponed to the
Appendix of the paper.

In subsection 
\ref{subs:basic} 
we write down an identity which turns out to be the stochastic
approximation of the PDE 
(\ref{eq:leroux}). 
Various error terms are defined here which will be
estimated in the forthcoming subsections.

In subsection 
\ref{subs:bounds} 
we introduce the relevant
\emph{Sobolev norms} 
and by using the previously proved a priori
estimates  we prove the necessary upper bounds on the apropriate
Sobolev norms of the error terms.

In subsection 
\ref{subs:mvsol} 
we show that choosing a subsequence of
the random Young measures 
(\ref{eq:young}) 
which converges vaguely in
distribution, the limit (random) Young measure is almost surely
measure valued entropy solution of the Cauchy problem
(\ref{eq:leroux}), 
(\ref{eq:ic}).

Subsection 
\ref{subs:coco} 
contains the stochastic version of the
method of 
\emph{compensated compactness}. 
It is further broken up into
two sub-subsections as follows. In sub-subsection
\ref{subsubs:murat}  
we preent the stochastic version of Murat's
Lemma: we prove that for any smooth Lax entropy/flux pair the entropy
production process is tight in the Sobolev space
$H^{-1}_{t,x}$.
In sub-subsection 
\ref{subsubs:tartar} 
we apply (an almost
sure version of) Tartar's Div-Curl Lemma leading to the desired
almost sure factorization property of the limiting random Young
measures. Finally, as main consequence of Tartar's Lemma,
we conclude that choosing any subsequence of
the random Young measures 
(\ref{eq:young}) 
which converges vaguely in
distribution, the limit (random) Young measure is almost surely
of Dirac type.

The  results of subsection 
\ref{subs:mvsol}  
and sub-subsection
\ref{subsubs:tartar} 
imply the Theorem. The
concluding steps are presented
in subsection 
\ref{subs:end}.

\subsection{Local equilibrium and a priori bounds}
\label{subs:apriori}

The hydrodynamic limit relies on macroscopically fast convergence to
(local) equilibrium in blocks of mesoscopic size 
$l$. 
Fix the block size  
$l$ 
and 
$(N,Z)\in\N\times\Z$
with the restriction 
$N+\abs{Z}\le l$ 
and denote  
\beqs
\Omega^l_{N,Z}
&:=&
\big\{
\uome\in\Omega^l\,:\,
\sum_{j=1}^l\eta_j=N, \sum_{j=1}^l\xi_j=Z
\big\},
\\[5pt]
\pil_{N,Z}(\underline\omega) 
&:=&
\pil_{\rho,u}(\underline\omega\,|\, \sum_{j=1}^l\eta_j=N,
\sum_{j=1}^l\xi_j=Z),
\eeqs
and, for 
$f:\Omega^l_{N,Z}\to\R$ 
\beqs
K^l_{N,Z}f(\uome)
&:=&
\sum_{j=1}^{l-1}
\big(f(\Theta_{j,j+1}\uome)-f(\uome)\big),
\\[5pt] 
D^l_{N,Z}(f)
&:=&
\frac12
\sum_{j=1}^{l-1}
\expect^l_{N,Z}
\left(
\big(f(\Theta_{j,j+1}\uome)-f(\uome)\big)^2
\right).
\eeqs
In plain words: 
$\Omega^l_{N,Z}$ 
is the hyperplane of configurations
$\uome\in\Omega^l$ 
with fixed values of the conserved quantities,
$\pi^l_{N,Z}$ 
is the  
\emph{microcanonical distribution} 
on this hyperplane,  
$K^l_{N,Z}$ 
is the symmetric infinitesimal generator
restricted to the hyperplane 
$\Omega^l_{N,Z}$,
and finally 
$D^l_{N,Z}$
is the Dirichlet form associated to 
$K^l_{N,Z}$. 
Note, that 
$K^l_{N,Z}$
is defined with free boundary conditions. 
Expectations with respect to the measures
$\pi^l_{N,Z}$ 
are denoted by 
$\expect^l_{N,Z}\big(\cdot\big)$. 
The convergence to local
equilibrium is 
\emph{quantitatively controlled}  
by the following uniform  logarithmic Sobolev estimate: 

\begin{lemma}
\label{lemma:lsi}
There exists a finite constant 
$\aleph$ 
such that for any 
$l\in\N$,
$(N,Z)\in\N\times\Z$ 
with 
$N+\abs{Z}\le l$ 
and any
$h:\Omega^l_{N,Z}\to\R_+$ 
with 
$\expect^l_{N,Z}(h)=1$ 
the following bound holds: 
\beq
\label{eq:lsi}
\expect^l_{N,Z}\big(h\log h\big)
\le
\aleph \, l^2
D^l_{N,Z} \left(\sqrt{h}\right).
\eeq
\end{lemma}

\noindent
{\bf Remark:}
In \cite{yau2} (see also \cite{leeyau}) the similar statement is
proved (inter alia) for symmetric simple exclusion process. That proof
can be easily adapted to our case. Instead of stirring  configurations
of two colours we have stirring of configurations of three colours. 
No really new ideas are involved. For sake of completeness however, we
sketch the proof in subsection \ref{subs:lsiproof} of the Appendix. 

\smallskip

The following  large deviation bound  goes back to Varadhan
\cite{varadhan}. See also the monographs \cite{kipnislandim} and
\cite{fritz1}. 

\begin{lemma}
\label{lemma:varadhan}
Let
$l\le n$, ${\cal V}:S^l\to \R_+$
and denote
${\cal V}_j(\uome) :=
{\cal V}(\omega_{j},\dots,\omega_{j+l-1})$.
Then for any
$\beta>0$
\beq
\label{eq:varadhan}
\frac{1}{n}  \sum_{j\in\Tn} \int_0^T
\expect_{\mun_s}\left( {\cal V}_j \right)\, ds
\le
C 
\frac{l^3}{\beta n^2\sigma}
+
\frac{T}{\beta}
\max_{ N, Z }
\log
\expect^{l}_{N,Z}
\big(\exp\left\{\beta {\cal V} \right\}\big)
\eeq
\end{lemma}

\noindent
{\bf Remarks:}
(1)  
Assuming only uniform bound of order $l^{-2}$ on the spectral gap of
$K^l_{N,Z}$ (rather than the stronger logarithmic Sobolev inequality
(\ref{eq:lsi})) and using Rayleigh-Schr\"odinger perturbation (see
Appendix 3 of \cite{kipnislandim}) we would get 
\beqs
&&
\hskip-17mm
\frac{1}{n}  \sum_{j\in\Tn} \int_0^T
\expect_{\mun_s}\left( {\cal V}_j \right)\, ds
\le
\\[5pt]
&&
C
\frac{l^3\Vert{\cal V}\Vert_\infty} {n^2\sigma}
+
T\Vert{\cal V}\Vert_\infty
\left(
\frac{\displaystyle\max_{N,Z} \expect^l_{N,Z}\big({\cal V}\big)}
     {\Vert{\cal V}\Vert_\infty}
+
\frac{\displaystyle\max_{N,Z} \var^l_{N,Z}\big({\cal V}\big)}
     {4\Vert{\cal V}\Vert_\infty^2}
\right),
\eeqs
which wouldn't be sufficient for our needs. 
\\
(2) 
The proof of the bound (\ref{eq:varadhan}) explicitly relies on the
logarithmic 
Sobolev inequality (\ref{eq:lsi}). It appears  in \cite{yau3}
and it is reproduced in several places, see e.g. \cite{fritz1,
  fritz2}. We do 
not repeat it here. 

\smallskip

The main probabilistic ingredients of our proof are  the following two 
consequences of Lemma 
\ref{lemma:varadhan}. 
These are variants of the celebrated 
\emph{one block estimate},
respectively, 
\emph{two blocks estimate} 
of Varadhan and  co-authors.

\begin{proposition}
\label{propo:apriori}
Assume conditions 
(\ref{cond:viscosity}) 
and  
(\ref{cond:blocksize}).
Given a local variable
$\ups_j$
there exists a constant
$C$
(depending
only on
$\ups_j$)
such that the following bounds hold:
\beq
\label{eq:obrepl}
\hskip-10mm
\expect
\Big(
\int_0^T\int_{\T}
\abs{
\wih{\ups}(s,x)
-
\Ups(\wih{\vxi}(s,x))
}^2
\,dx\,dt
\Big)
&\le&
C\frac{l^2}{n^2\sigma}
\\[5pt]
\label{eq:gradbound}
\hskip-10mm
\expect
\Big(
\int_0^T\int_{\T}
\abs{
\px \wih{\ups}(s,x)
}^2
\,dx\,dt
\Big)
&\le&
C
\sigma^{-1}
\eeq
\end{proposition}

The proof of Proposition \ref{propo:apriori} is postponed to
subsection \ref{subs:aprioriproof} in the Appendix. It relies on the
large  deviation bound (\ref{eq:varadhan}) and an  elementary
probability  lemma stated in subsection \ref{subs:epl} of the
Appendix.

We shall refer to  
(\ref{eq:obrepl}) 
as the 
\emph{block replacement  bound} 
and to 
(\ref{eq:gradbound}) 
as the 
\emph{gradient  bound}.

\subsection{The basic identity}
\label{subs:basic}

Given a smooth function
$f:\dom \to \R$
we write
\beqs
\pt f(\wih{\vxi}(t,x))
=
\Gn f(\wih{\vxi}(t,x))
+
\pt M^n_f(t,x),
\eeqs
where the process
$t\mapsto M^n_f(t,x)$
is a martingale. Here and in
the future
$\pt f(\wih{\vxi}(t,x))$
and
$\pt M^n_f(t,x)$
are meant as
\emph{distributions} 
in their time variable.

In this order we compute the action of the
infinitesimal generator
$\Gn=n\Ln + n^2\sigma\Kn$
on
$f(\wih{\vxi}(x))$.
First we compute the asymmetric part:
\beq
\label{eq:lf}
n\Ln f(\wih{\vxi}(x))
&=&
-
\nabla f(\wih{\vxi}(x))\cdot\px\wih{\vphi}(x)
+
A^{1,n}_f(x)
\eeq
where
\beq
\label{eq:a1}
A^{1,n}_f(x)
=
A^{1,n}_f(\uome,x)
:=
n\sum_{j\in\T}
r(\omega_j,\omega_{j+1})\times
\hskip4.5cm
\\[5pt]
\notag
\Big\{
f\big(
\wih{\vxi}(x)
-
\frac{1}{l}
\big(
a(\frac{nx-j}{l})-a(\frac{nx-j-1}{l})
\big)
\big( \vxi_j-\vxi_{j+1}\big)
\big)
-
f\big(\wih{\vxi}(x)\big)
\phantom{\Big\}}
\\[5pt]
\notag
\phantom{.}
\hskip3cm
+
\frac{1}{l^2} \, a'(\frac{nx-j}{l})
\nabla f\big(\wih{\vxi}(x)\big)
\cdot\big(\vxi_j-\vxi_{j+1}\big)
\Big\}.
\eeq
See formula 
(\ref{eq:phipsidef}) 
for the definition of
$\vphi$.
$A^{1,n}_f$ 
is a numerical  error term which will be easy to estimate.

Next, the symmetric part:
\beq
\label{eq:kf}
n^2\sigma\Kn f(\wih{\vxi}(x))
&=&
\sigma\nabla f(\wih{\vxi}(x))\cdot\px^2\,\wih{\vxi}(x)
+
A^{2,n}_f(x)
\eeq
where
\beq
\label{eq:a2}
A^{2,n}_f(x)
=
A^{2,n}_f(\uome,x)
:=
\hskip7cm
\\[5pt]
\notag
n^2\sigma
\sum_{j\in\T}
\Big\{
f\big(
\wih{\vxi}(x)
-
\frac{1}{l}
\big(
a(\frac{nx-j}{l})-a(\frac{nx-j-1}{l})
\big)
\big( \vxi_j-\vxi_{j+1}\big)
\big)
\phantom{\Big\}}
\\[5pt]
\notag
\phantom{.}
\hskip3cm
-
f\big(\wih{\vxi}(x)\big)
\phantom{\Big\}}
+
\frac{1}{l^3} a''(\frac{nx-j}{l})
\nabla f\big(\wih{\vxi}(x)\big)
\cdot\vxi_j
\Big\}.
\eeq
This is another numerical error term easy to estimate.

Hence our 
\emph{basic identity}
\beq
\label{eq:basic}
&&
\hskip-5mm
\pt f(\wih{\vxi}(t,x))
+
\nabla f(\wih{\vxi}(t,x))
\cdot
\nabla \vPhi (\wih{\vxi}(t,x))
\cdot
\px\wih{\vxi}(t,x)
=
\\[5pt]
&&
\notag
\hskip5mm
\sum_{i=1}^2
\left(
A^{i,n}_f(t,x)
+
B^{i,n}_f(t,x)
+
C^{i,n}_f(t,x)
\right)
+
\pt M^{n}_f(t,x)\,.
\eeq
The various terms on the right hand side are
\beq
\label{eq:b1}
&&
\hskip-8mm
B^{1,n}_f(x)
=
B^{1,n}_f(\uome,x)
:=
\px
\left\{
\nabla f(\wih{\vxi}(x))
\cdot
\big(
\vPhi(\wih{\vxi}(x))-\wih{\vphi}(x)
\big)
\right\}
\\[5pt]
\label{eq:b2}
&&
\hskip-8mm
B^{2,n}_f(x)
=
B^{2,n}_f(\uome,x)
:=
\sigma
\px^2f(\wih{\vxi}(x))
=
\px
\left\{
\sigma
\nabla f(\wih{\vxi}(x))
\cdot
\px \wih{\vxi}(x)
\right\}
\\[5pt]
\notag
&&
\hskip-8mm
C^{1,n}_f(x)
=
C^{1,n}_f(\uome,x)
:=
-
\big(\px\wih{\vxi}(x)\big)^\dagger
\cdot
\nabla^2f(\wih{\vxi}(x))
\cdot
\big(
\vPhi(\wih{\vxi}(x))-\wih{\vphi}(x)
\big)
\\
\label{eq:c1}
\\
\label{eq:c2}
&&
\hskip-8mm
C^{2,n}_f(x)
=
C^{2,n}_f(\uome,x)
:=
-
\sigma
\big(\px\wih{\vxi}(x)\big)^\dagger
\cdot
\nabla^2f(\wih{\vxi}(x))
\cdot
\big(\px\wih{\vxi}(x)\big)
\eeq
and
\beqs
&&
A^{i,n}_f(t,x):=A^{i,n}_f(\Xn_t,x),
\\[5pt]
&&
B^{i,n}_f(t,x):=B^{i,n}_f(\Xn_t,x),
\\[5pt]
&&
C^{i,n}_f(t,x):=C^{i,n}_f(\Xn_t,x).
\eeqs

In the present paper we shall apply the basic identity
(\ref{eq:basic}) 
only for Lax entropies
$f(\vu)=S(\vu)$.
In this special case the left hand side gets the form of a
conservation law: 
\beq
\label{eq:lentbasic}
&&
\hskip-5mm
\pt S(\wih{\vxi}(t,x))
+
\px F(\wih{\vxi}(t,x))
=
\\[5pt]
&&
\notag
\hskip5mm
\sum_{i=1}^2
\left(
A^{i,n}_S(t,x)
+
B^{i,n}_S(t,x)
+
C^{i,n}_S(t,x)
\right)
+
\pt M^{n}_S(t,x),
\eeq

\subsection{Bounds}
\label{subs:bounds}

We fix
$T<\infty$
and use the
$L^p$
norms
\beqs
\norm{g}_{L^{p}_{t,x}}^p
:=
\int_0^T\int_{\T}
\abs{g(t,x)}^p\,dx\,dt
\eeqs
and the Sobolev norms
\[
\norm{g}_{W^{-1,p}_{t,x}}
:=
\sup\big\{
\int_0^T \int_{\T}
\varphi(t,x) g(t,x) \,dx\,dt\,:\,
\norm{\pt\varphi}_{L^{q}_{t,x}}^q +
\norm{\px\varphi}_{L^{q}_{t,x}}^q \le1 
\big\}
\]
where
$p^{-1}+q^{-1}=1$
and
$\varphi:[0,T]\times\T\to \R$
is a test function.  We use   the standard  notation
$W_{t,x}^{-1,2}=:H_{t,x}^{-1}$.

\smallskip
\noindent
{\bf Remark on notation:}
The numerical error terms
$A^{i,n}_f(t,x)$, $i=1,2$,
will be estimated in
$L^{\infty}_{t,x}$
norm. In these estimates only Taylor expansion bounds are used, no
probabilistic argument is involved. The more sophisticated terms
$B^{i,n}_f(t,x)$,  $i=1,2$,
respectively,
$C^{i,n}_f(t,x)$, $i=1,2$,
will be estimated in
$H^{-1}_{t,x}$,
respectively,
$L^{1}_{t,x}$
norms. The martingale derivative
$\pt M^n_f(t,x)$
will be estimated in
$H^{-1}_{t,x}$
norm.

\smallskip

By straightforward numerical estimates (which do not rely on any
probabilistic arguments) we obtain

\begin{lemma}
\label{lemma:bounds1}
Assume conditions 
(\ref{cond:viscosity}) 
and 
(\ref{cond:blocksize}).
Let
$f:\dom\to\R$
be a twice continuously differentiable function with
bounded derivatives. Then almost surely
\beqs
\norm{A^{1,n}_f}_{L^{\infty}_{t,x}}
=
o(1)\quad\text{and}\quad
\norm{A^{2,n}_f}_{L^{\infty}_{t,x}}
=
o(1)
\eeqs
as $n\to\infty$.
\end{lemma}

\begin{proof}
Indeed, using nothing more than Taylor expansion and boundedness of
the local variables we readily obtain
\beq
\label{eq:a1bound}
&&
\sup_{x\in\T}\sup_{\uome\in\Omn}
\abs{A^{1,n}_f(\uome,x)}
\le
C\frac{n}{l^{2}}
=o(1)
\\[5pt]
\label{eq:a2bound}
&&
\sup_{x\in\T}\sup_{\uome\in\Omn}
\abs{A^{2,n}_f(\uome,x)}
\le
C\frac{n^2\sigma}{ l^{3}}
=o(1).
\eeq
We omit the tedious but otherwise straightforward details.
\end{proof}

Applying Proposition
\ref{propo:apriori} 
we obtain the following more sophisticated bounds

\begin{lemma}
\label{lemma:bounds2}
Assume conditions 
(\ref{cond:viscosity}) 
and  
(\ref{cond:blocksize}).
Let
$f:\dom\to\R$
be a twice continuously differentiable function with
bounded derivatives. The following asymptotics hold, as
$n\to\infty$:
\beqs
\phantom{.}
(i)
&\phantom{MMMMMMMM}&
\expect\left(\norm{B^{1,n}_f}_{H^{-1}_{t,x}} \right)
=o(1)
\phantom{MMMMMMMMM}
\\[5pt]
(ii)
&&
\expect\left(\norm{B^{2,n}_f}_{H^{-1}_{t,x}} \right)
=o(1)
\\[5pt]
(iii)
&&
\expect\left(\norm{C^{1,n}_f}_{L^{1}_{t,x}} \right)
=o(1)
\\[5pt]
(iv)
&&
\expect\left(\norm{C^{2,n}_f}_{L^{1}_{t,x}} \right)
=\Ordo(1)
\eeqs
\end{lemma}

\begin{proof}
\phantom{.} \hfill
\\
\noindent
$(i)$
We use the block replacement bound 
(\ref{eq:obrepl}):
\beqs
&&
\hskip-7mm
\expect
\Big(
\abs{
\int_0^T\int_{\T}
v(t,x)B^{1,n}_f(t,x)
\,dx\,dt
}
\Big)
\\[5pt]
&&
=
\expect
\Big(
\abs{
\int_0^T\int_{\T}
\px v(t,x)
\nabla f(\wih{\vxi}(t,x))
\cdot
\big(
\vPhi(\wih{\vxi}(t,x))-\wih{\vphi}(t,x)
\big)
\,dx\,dt
}
\Big)
\\[5pt]
&&
\le
\sup_{\vu\in{\cal D}}\abs{\nabla f(\vu)}
\norm{\px v}_{L^2_{t,x}}
\expect
\Big(
\int_0^T\int_{\T}
\abs{
\vPhi(\wih{\vxi}(t,x))-\wih{\vphi}(t,x)
}^2
\,dx\,dt
\Big)^{1/2}
\\[5pt]
&&
\le
C
\norm{\px v}_{L^2_{t,x}}
\frac{l}{n\sqrt{\sigma}}.
\eeqs
\noindent
$(ii)$
We use the gradient bound 
(\ref{eq:gradbound}):
\beqs
&&
\hskip-12mm
\expect
\Big(
\abs{
\int_0^T\int_{\T}
v(t,x)B^{2,n}_f(t,x)
\,dx\,dt
}
\Big)
\\[5pt]
&&
=
\expect
\Big(
\abs{
\int_0^T\int_{\T}
\px v(t,x)
\nabla f(\wih{\vxi}(t,x))
\cdot
\sigma
\big(
\px\wih{\vxi}(t,x)
\big)
\,dx\,dt
}
\Big)
\\[5pt]
&&
\le
\sup_{\vu\in{\cal D}}\abs{\nabla f(\vu)}
\norm{\px v}_{L^2_{t,x}}
\sigma
\expect
\Big(
\int_0^T\int_{\T}
\abs{
\px\wih{\vxi}(t,x)
}^2
\,dx\,dt
\Big)^{1/2}
\\[5pt]
&&
\le
C
\norm{\px v}_{L^2_{t,x}}
\sigma^{1/2}.
\eeqs
\noindent
$(iii)$
We use both, the block replacement bound
(\ref{eq:obrepl}) 
and the gradient bound 
(\ref{eq:gradbound}):
\beqs
&&
\hskip-10mm
\expect
\Big(
\int_0^T\int_{\T}
\abs{
C^{1,n}_f(t,x)
}
\,dx\,dt
\Big)
\\[5pt]
&&
\le
\sup_{\vu\in{\cal D}}\abs{\nabla^2f(\vu)}
\expect
\Big(
\int_0^T\int_{\T}
\abs{
\wih{\vphi}(s,x)
-
\vPhi(\wih{\vxi}(s,x))
}^2
\,dx\,dt
\Big)^{1/2}
\times
\\[5pt]
&&
\hskip28mm
\expect
\Big(
\int_0^T\int_{\T}
\abs{
\px \wih{\vxi}(s,x)
}^2
\,dx\,dt
\Big)^{1/2}
\\[5pt]
&&
\le
C
\frac{l}{n\sigma}.
\eeqs
\noindent
$(iv)$
We use again the gradient bound 
(\ref{eq:gradbound}):
\beqs
&&
\hskip-18mm
\expect
\Big(
\int_0^T\int_{\T}
\abs{
C^{2,n}_f(t,x)
}
\,dx\,dt
\Big)
\\[5pt]
&&
\le
\sup_{\vu\in{\cal D}}\abs{\nabla^2f(\vu)}
\sigma
\expect
\int_0^T\int_{\T}
\abs{
\px \wih{\vxi}(s,x)
}^2
\,dx\,dt
\\[5pt]
&&
\le
C.
\eeqs
\end{proof}

\begin{lemma}
\label{lemma:bounds3}
Assume conditions 
(\ref{cond:viscosity}) 
and  
(\ref{cond:blocksize}).
Let
$f:\dom\to\R$
be a twice continuously differentiable function with
bounded derivatives. There exists a constant
$C$
(depending only on
$f$)
such that the folowing asymptotics holds as
$n\to\infty$:
\beqs
\expect\left(\norm{\pt M^{n}_f}_{H^{-1}_{t,x}} \right)
=o(1)
\eeqs
\end{lemma}

\begin{proof}
Since
\beqs
\norm{\pt M^{n}_f}^2_{H^{-1}_{t,x}}
\le
\norm{M^{n}_f}^2_{L^{2}_{t,x}},
\eeqs
we have to bound  the expectation of the right hand side.
\beqs
\expect
\Big(
\int_0^T\int_\T
\big( M^{n}_f(t,x) \big)^2
\,dx\,dt
\Big)
=
\expect
\Big(
\int_0^T\int_\T
\langle M^{n}_f(t,x) \rangle
\,dx\,dt
\Big),
\eeqs
where
$t\mapsto \langle M^{n}_f(t,x) \rangle$
is the conditional
variance process of the martingale
$M^{n}_f(t,x)$:
\beqs
\label{eq:condcov}
\langle
M^{n}_f(t,x)
\rangle
&=&
\phantom{++}\,
n
\left(
\Ln f^2(\wih{\vxi}(t,x))
-
2f(\wih{\vxi}(t,x))\Ln f(\wih{\vxi}(t,x))
\right)
\\[5pt]
&&
+
n^2\sigma
\left(
\Kn f^2(\wih{\vxi}(t,x))
-
2f(\wih{\vxi}(t,x))\Kn f(\wih{\vxi}(t,x))
\right).
\eeqs
Using the expressions 
(\ref{eq:lf}) 
and 
(\ref{eq:kf}) 
we obtain
\beqs
\langle
M^{n}_f(t,x)
\rangle
&=&
\phantom{+}
A^{1,n}_{f^2}(t,x)-2f(\wih{\vxi}(t,x))A^{1,n}_f(t,x)
\\
&&
+
A^{2,n}_{f^2}(t,x)-2f(\wih{\vxi}(t,x))A^{2,n}_f(t,x).
\eeqs
Hence, by the bounds 
(\ref{eq:a1bound}) 
and 
(\ref{eq:a2bound}) 
(which apply as well of course to the function
$f^2$),
we obtain
\[
\sup_{t\in[0,T]}
\sup_{x\in\T}
\,
\langle
M^{n}_f(t,x)
\rangle
\le
C\frac{n^2\sigma}{l^3}
=o(1),
\]
which proves the lemma.
\end{proof}

\subsection{Convergence to measure valued entropy solutions}
\label{subs:mvsol}

\begin{proposition}
\label{prop:mvsol}
Conditions 
(\ref{cond:viscosity}), 
(\ref{cond:blocksize}),
and 
(\ref{cond:iniprofile}) 
are in force. Let
$\Q^{n'}$
be a subsequence of the
probability distributions
defined in 
(\ref{eq:young}), 
which converges weakly in the vague
sense:
$\Q^{n'}\wvto\Q$.
Then the probability measure
$\Q$
is concentrated on the  measure valued entropy solutions of the Cauchy
problem   
(\ref{eq:leroux}), 
(\ref{eq:ic}).
\end{proposition}

\begin{proof}
Due to separability of
$C([0,T]\times\T)$
it is sufficient to prove
that for any convex  Lax entropy/flux pair
$(S,F)$
and any positive test function
$\varphi:[0,T]\times\T\to\R_+$,
(\ref{eq:mlax})
holds
$\Q$-almost-surely.
So we fix
$(S,F)$
and
$\varphi$,
and denote the real random variable
\beqs
X^n
&:=&
-
\int_0^T\int_\T
\varphi(t,x)
\big(
\pt S(\wih{\vxi}(t,x))
+
\px F(\wih{\vxi}(t,x))
\big)
\,dx\,dt
\\[5pt]
&=&
\phantom{-}
\int_0^T\int_\T
\big(
\pt\varphi(t,x)
\langle\nu^n(t,x)\,,\, S\rangle 
+ 
\px\varphi(t,x)
\langle\nu^n(t,x)\,,\, F\rangle 
\big)
\,dx\,dt
\\[5pt]
&&
\hskip3.9cm
+
\int_\T
\varphi(0,x)
\langle\nu^n(0,x)\,,\, S\rangle
\,dx.
\eeqs
The right hand side is a continuous function of
$\nu^n$,
so from the asumption
$\Q^n\wvto\Q$
it follows that
\beq
\label{eq:xtox}
X^n\Rightarrow X,
\eeq
where
\beqs
X
&:=&
\phantom{-}
\int_0^T\int_\T
\big(
\pt\varphi(t,x)
\langle\nu(t,x)\,,\,S\rangle 
+ 
\px\varphi(t,x)
\langle\nu(t,x)\,,\,F\rangle 
\big)
\,dx\,dt
\\[5pt]
&&
\hskip3.9cm
+
\int_\T
\int_{\cal D}
\nu(0,x;d\vv)
\varphi(0,x)S(\vv)
\,dx.
\eeqs
and
$\nu$
is distributed according to
$\Q$.

We apply the basic identity 
(\ref{eq:basic}) 
speciafied for
$f(\vu)=S(\vu)$, 
that is identity (\ref{eq:lentbasic}). It follows that
\beq
\label{eq:decompo}
X^n=
Y^n+Z^n
\eeq
where
\beqs
Y^n
&:=&
\int_0^T\int_\T
\varphi(t,x)
C^{2,n}_S(t,x)
\,dx\,dt
\\[5pt]
&=&
\notag
\sigma
\int_0^T\int_\T
\varphi(t,x)
\big(\px\wih{\vxi}(t,x)\big)^\dagger
\cdot
\nabla^2 S(\wih{\vxi}(t,x))
\cdot
\big(\px\wih{\vxi}(t,x)\big)
\eeqs
and
\beqs
Z^n
:=
\int_0^T\int_\T
\varphi(t,x)
\left(
\sum_{i=1}^2
\left(
A^{i,n}_S
+
B^{i,n}_S
\right)
+
C^{1,n}_S
+
\pt M^n_S
\right)
(t,x)
\,dx\,dt.
\eeqs
Due to convexity of
$S$
and positivity of
$\varphi$
we have
\beq
\label{eq:ypositive}
Y^n\ge0,
\qquad
\text{ almost surely.}
\eeq
On the other hand, from Lemmas
\ref{lemma:bounds1}, 
\ref{lemma:bounds2}, 
\ref{lemma:bounds3}
we conclude that
\beq
\label{eq:ztozero}
\lim_{n\to\infty}
\expect\big(\abs{Z^n}\big)
=0.
\eeq
Finally, from 
(\ref{eq:xtox}), 
(\ref{eq:decompo}),
(\ref{eq:ypositive}) 
and 
(\ref{eq:ztozero}) 
the statement of the Proposition follows.
\end{proof}

\subsection{Compensated compactness}
\label{subs:coco}

\subsubsection{Murat's lemma}
\label{subsubs:murat}

\begin{lemma}
\label{lemma:murat}
Assume conditions 
(\ref{cond:viscosity}) 
and 
(\ref{cond:blocksize}).
Given a twice continuously differentiable Lax entropy/flux pair
$(S,F)$,
the sequence
\beqs
 X^{n}(t,x):=
\pt S(\wh\vxi^{n}(t,x))
+
\px F(\wh\vxi^{n}(t,x))
\eeqs
is tight in
$H^{-1}_{t,x}$.
\end{lemma}

\begin{proof}

Note that 
$ X^{n}(t,x)$
is exactly 
the left hand side of the basic identity 
(\ref{eq:lentbasic})
and racall that this expression (in particular
$\pt S(\wih{\vxi}(t,x))$)
is a  random distribution in its
$t$
variable.

By definition and a priori boundedness of the domain
$\dom$,
there exists a constant
$C<\infty$
such that
\beq
\label{eq:murat1}
\prob
\Big(
\norm{X^{n}}_{W^{-1,\infty}_{t,x}}\le C
\Big)
=1.
\eeq
We decompose
\beq
\label{eq:decomp}
X^{n}(t,x)
=
Y^{n}(t,x)+Z^{n}(t,x),
\eeq
where
\beqs
&&
Y^{n}(t,x)
:=
B^{1,n}_S(t,x)
+
B^{2,n}_S(t,x)
+
\pt M^{n}_S(t,x),
\\[5pt]
&&
Z^{n}(t,x)
:=
A^{1,n}_S(t,x)
+
A^{2,n}_S(t,x)
+
C^{1,n}_S(t,x)
+
C^{2,n}_S(t,x).
\eeqs
For the definitions of the terms 
$A^{i,n}_S$, 
$B^{i,n}_S$, 
$C^{i,n}_S$, 
$i=1,2$, 
see 
(\ref{eq:a1}), 
(\ref{eq:a2})
and 
(\ref{eq:b1})--(\ref{eq:c2}).
   
From Lemmas 
\ref{lemma:bounds1},
\ref{lemma:bounds2} 
and 
\ref{lemma:bounds3} 
it follows that
\beq
\label{eq:murat2}
\expect
\Big(
\norm{Y^{n}}_{H^{-1}_{t,x}}
\Big)
\to0,
\eeq
and
\beq
\label{eq:murat3}
\expect
\Big(
\norm{Z^{n}}_{L^{1}_{t,x}}
\Big)
\le C.
\eeq
Further on, from 
(\ref{eq:murat2}), 
respectively, 
(\ref{eq:murat3})
it follows that for any
$\vareps>0$
one can find a 
\emph{compact} 
subset
$K_\vareps$ of $H^{-1}_{t,x}$
and a 
\emph{bounded} 
subset
$L_\vareps$ 
of 
$L^1_{t,x}$
such that
\beq
\label{eq:smallprob}
\prob\Big(Y^n\notin K_\vareps\Big)<\vareps/2,
\qquad
\prob\Big(Z^n\notin L_\vareps\Big)<\vareps/2.
\eeq
On the other hand,  Murat's lemma (see 
\cite{murat} 
or Chapter 9 of
\cite{serre})  
says that
\[
M_\vareps:=
\big(K_\vareps + L_\vareps\big)\cap
\{X\in H^{-1}_{t,x}: \norm{X}_{ W^{-1,\infty}_{t,x}}\le C\}
\]
is compact in
$H_{t,x}^{-1}$.
From 
(\ref{eq:murat1}), 
(\ref{eq:decomp}) 
and
(\ref{eq:smallprob}) 
it follows that
\[
\prob\Big(X^n\notin M_\vareps\Big)<\vareps,
\]
uniformly in
$n$,
which proves the lemma.
\end{proof}

\subsubsection{Tartar's lemma and its consequence}
\label{subsubs:tartar}

\begin{lemma}
\label{lemma:tartar}
Assume conditions 
(\ref{cond:viscosity}) 
and  
(\ref{cond:blocksize}).
Let
$\Q^{n'}$
be a subsequence of the probability measures on
$\cal Y$
defined in 
(\ref{eq:young}), 
which converges weakly in the vague
sense:
$\Q^{n'}\wvto\Q$.
Then
$\Q$
is concentrated on the (vaguely closed)  subset of Young measures
satisfying 
\eqref{tarfac}.
That is, $\Q$-a.s.  
for any two generalized Lax entropy/flux pairs
$(S_1,F_1)$
and
$(S_2,F_2)$
and any test function
$\varphi:[0,T]\times\T\to\R$,
\beq
\label{eq:tartar}
&&
\hskip-6mm
\int_0^T\int_\T
\varphi(t,x)
\langle \nu(t,x) \,,\, S_1F_2-S_2F_1 \rangle
\,dx\,dt
=
\\[5pt]
\notag
&&
\hskip-3mm
\int_0^T\int_\T
\varphi(t,x)
\big(
\langle \nu(t,x) \,,\, S_1\rangle
\langle \nu(t,x) \,,\, F_2\rangle
-
\langle \nu(t,x) \,,\, S_2\rangle
\langle \nu(t,x) \,,\, F_1\rangle
\big)
\,dx\,dt.
\eeq
\end{lemma}

\begin{proof}
First we prove (\ref{eq:tartar}) for twice continuously
differentiable entropy/flux pairs. 
Due to separability of
$C([0,T]\times\T)$
it is sufficient to prove that for any two twice continuously
differentiable Lax entropy/flux pairs
$(S_1,F_1)$
and
$(S_2,F_2)$
and any test function
$\varphi:[0,T]\times\T\to\R$,
(\ref{eq:tartar}) 
holds
$\Q$-almost-surely.
So we fix
$(S_1,F_1)$, 
$(S_2,F_2)$
and
$\varphi$.
Note that
\beqs
X^n_j(t,x)
&:=&
\pt S_j(\wih{\vxi}(t,x))
+
\px F_j(\wih{\vxi}(t,x))
\\[5pt]
&=&
\pt
\langle\nu^n(t,x)\,,\,S_j\rangle
+
\px
\langle\nu^n(t,x)\,,\,F_j\rangle
\eeqs
$j=1,2$.

Due to Skorohod's representation theorem (see Theorem 1.8 of
\cite{ethierkurtz})
and Lemma \ref{lemma:murat} 
we can realize the random Young measures
$\nu^n(t,x;d\vv)$
and
$\nu(t,x;d\vv)$
\emph{jointly} 
on an enlarged probablity space
$(\Xi,{\cal A}, \prob)$
so that 
\emph{$\prob$-almost-surely}
\beqs
\nu^{n'} \vto \nu,
\qquad
\text{ and }
\qquad
\{X^{n'}_j:\,{n'},\,\, j=1,2\, \} 
\text{ is relative compact in }
H^{-1}_{t,x}.
\eeqs
So, applying Tartar's Div-Curl Lemma (see 
\cite{tartar1},
\cite{tartar2}, 
or Chapter 9 of 
\cite{serre}) 
we conclude that (in
this realization) almost surely the factorization 
(\ref{eq:tartar})
holds true. 

Since $\dom$ is compact, from Riemann's method of solving the
linear hyperbolic PDE (\ref{entropia2}) (see Chapter 4 of
\cite{john}) it follows that
generalized entropy/flux pairs are 
approximated pointwise by
smooth ones. Thus the Tartar factorization  (\ref{eq:tartar})
extends from smooth to generalized entropy/flux pairs. 
Hence the lemma.
\end{proof}

The main consequence of Lemma \ref{lemma:tartar} is the following

\begin{proposition}
\label{prop:diperna}
Assume conditions 
(\ref{cond:viscosity}) 
and  
(\ref{cond:blocksize}).
Let
$\Q^{n'}$
be a subsequence of the probability measures on
$\cal Y$
defined in 
(\ref{eq:young}), 
which converges weakly in the vague sense:
$\Q^{n'}\wvto\Q$.
Then the probability measure
$\Q$
is concentrated on a set of Dirac-type Young measures, that is
$\Q({\cal U})=1$.
\end{proposition}

\begin{proof}
In view of Lemma 
\ref{lemma:tartar} 
this is a direct consequence of Lemma 
\ref{lemma:dirac}.
\end{proof}

\noindent{\bf Remark:}
This is the only point where we exploit the very special features of
the PDE  
(\ref{eq:leroux}). 
Note that the proof of Lemma \ref{lemma:dirac} relies on elementary
explicit  computations. In case of general
$2\times2$
hyperbolic systems of  conservation laws, instead of these explicit 
computations we should  refer to DiPerna's arguments from 
\cite{diperna}, 
possibly further complicated  by the existence of  singular
(non-hyperbolic)  points isolated 
at the boundary of the domain
$\cal D$.
More general results will be presented in the forthcoming paper
\cite{fritztoth2}.

\subsection{End of proof}
\label{subs:end}

From Propositions 
\ref{prop:mvsol} 
and 
\ref{prop:diperna} 
it follows that from any subsequence
$n'$
one can extract a sub-subsequence
$n''$
such that
$\Q^{n''}\wvto\Q$
and
$\Q$
is concentrated on the
set of Dirac-type measure valued entropy solutions of the Cauchy
problem. From now on we denote simply by
$n$
this sub-subsequence.
Referring again to Skorohod's Representation Theorem we
realize the Dirac-type random Young measures
$\nu^{n}_{t,x}(d\vv):=\delta_{\wih{\vxi}(t,x)}(d\vv)$
and
$\nu_{t,x}(d\vv):=\delta_{\vu(t,x)}(d\vv)$
jointly on an enlarged probability space
$(\Xi,{\cal A}, \prob)$,
so that
$\nu^{n}\vto\nu$
almost surely and
$(t,x)\mapsto\vu(t,x)$
is almost surely entropy solution of the
Cauchy problem.
From basic functional analytic
considerations (see e.g. Chapter 9 of 
\cite{serre}) 
it follows that, in case that the limit Youg measure is also
Dirac-type, the vague convergence
$\nu^{n}\vto\nu$
implies  strong (i.e. norm) convergence of the underlying
functions,
\beq
\label{eq:l1cvg}
\wih{\vxi}\to \vu
\quad\text{ in }\quad
L^1_{t,x}.
\eeq
So, we have realized jointly on the probability space
$(\Xi,{\cal A}, \prob)$
the empirical block average processes
$\wih{\vxi}(t,x)$
and the random function
$\vu(t,x)$
so that  the latter one
is almost surely entropy solution of the Cauchy problem, and
(\ref{eq:l1cvg}) 
almost surely holds true. This proves the theorem.
\hfill\hfill\qed

\section{Appendix}
\label{section:app}

\subsection{The logarithmic Sobolev inequality for random stirring of $r$
colours on the linear graph $\{1,2,\dots,l\}$}
\label{subs:lsiproof}

Let  $r\ge2$ be a fixed intger. For  $l\in\N$ we consider $r$-tuples
of  integers 
$N=(N_1,\dots,N_r)$ such that  
\beq
\label{eq:ncond}
&&
N_\alpha\ge0, \quad \alpha=1,\dots,r
\qquad\text{ and }\qquad
N_1+\dots+N_r=l,
\\[5pt]
\notag
&&
\Omega^l_N
:=
\big\{
\uome\in\{1,\dots,r\}^l:
\sum_{j=1}^l\ind_{\{\ome_j=\alpha\}}=N_\alpha, \,
\alpha=1,\dots,r
\big\}.
\eeq
Let $\pi^l_N$ denote the uniform probability measure on $\Omega^l_N$:
\[
\pi^l_N(\uome)=
\frac{N_1!\cdots N_r!}{l!}, 
\qquad
\uome\in\Omega^l_N. 
\]
The one  dimensional marginals of $\pi^l_N$ are
\[
\pi^{l,1}_N(\alpha)
=
\frac{N_\alpha}{l}.
\]
The random element of $\Omega^l_N$ distributed according to  $\pi^l_N$
will be denoted $\uzet=(\zeta_1,\zeta_2,\dots,\zeta_l)$. Expectation
with respect to $\pi^l_N$, respectively, $\pi^{l,1}_N$ 
will be denoted by
$\expect^l_N\big(\cdots\big)$, 
respectively, 
$\expect^{l,1}_N\big(\cdots\big)$. Conditional
expectation, given the first coordinate $\zeta_1$ will be
denoted 
$\expect^l_N\big(\cdots\big|\zeta_1\big)$. 
Note that 
\[
\expect^{l}_{N}\big(f(\uzet)\big|\zeta_1=\alpha \big)
=
\expect^{l-1}_{N^{\alpha}}\big(f(\alpha,\zeta_2,\dots,\zeta_l)\big)
\]
where 
$\expect^{l-1}_{N^{\alpha}}\big(\cdots\big)$ stands for
expectation with respect to $(\zeta_2,\dots,\zeta_l)$ distributed
according to $\pi^{l-1}_{N^{\alpha}}$ and, given
$N=(N_1,\dots,N_\alpha,\dots,N_r)$ with $N_\alpha\ge1$, $N^\alpha:=
(N_1,\dots,N_\alpha-1,\dots,N_r)$.

Given a probability density $h$ over $(\Omega^l_N,\pi^l_N)$, its
entropy is 
\beqs
H^l_N\big(h\big)
:=
\expect^l_N\big(h(\uzet)\log h(\uzet)\big).
\eeqs

Further on, for $i,j\in\{1,\dots,l\}$ let
$\Theta_{i,j}:\Omega^l_N\to\Omega^l_N$ be the spin exchange operator
\[
\left(\Theta_{i,j}\uome\right)_k=
\left\{
\begin{array}{ll}
\omega_j\quad&\text{ if } k=i, 
\\
\omega_i\quad&\text{ if } k=j, 
\\
\omega_k\quad&\text{ if } k\not=i,j,
\end{array}
\right..
\]
For $f:\Omega^l_N\to\R$ we define the Dirichlet form and the
conditional Dirichlet form, given $\zeta_1$
\beqs
D^l_N\big(f\big)
&:=&
\frac12
\sum_{i=1}^{l-1}
\expect^l_N
\left(
\big(f(\Theta_{i,i+1}\uzet)-f(\uzet)\big)^2
\right),
\\[5pt]
D^l_N\big(f\big|\zeta_1\big)
&:=&
\frac12
\sum_{i=1}^{l-1}
\expect^l_N
\big(
(f(\Theta_{i,i+1}\uzet)-f(\uzet))^2
\big|
\zeta_1
\big)
\\[5pt]
&=&
D^{l-1}_{N^{\zeta_1}}\big(f(\zeta_1,\cdot)\big).
\eeqs

The logarithmic Sobolev inequality is formulated in the following 

\begin{proposition}
There exist a finite constant 
$\aleph$ 
such that for any number of colours
$r$,
any block size 
$l\in \N$, 
any distribution of colours 
$N=(N_1,\dots,N_r)$
satisfying 
(\ref{eq:ncond})
and any probability density 
$h$ 
over 
$(\Omega^l_N, \pi^l_N)$, 
the following inequality holds: 
\beq
\label{eq:lsigen}
H^l_N\big(h\big)
\le 
\aleph\, l^2 
D^l_N\big(\sqrt h\big).
\eeq
\end{proposition}

\noindent
{\bf Remark:} 
The proof follows  \cite{yau2} (see also
\cite{leeyau}). Due to exchangeability of the measures $\pi^l_N$ 
some steps are considearbly simpler than there. 

\begin{proof}
We shall prove the Proposition by induction on $l$. 
Denote 
\[
W(l):=
\sup_{N}\sup_{h}
\frac{H^l_N\big(h\big)}{D^l_N\big(\sqrt h\big)}. 
\]
The following identity is straightforward
\beq
\notag
H^l_N\big(h\big)
&=&
\phantom{+}
\expect^{l,1}_N
\big(
\expect^l_N\big(h(\uzet) \big| \zeta_1 \big)
\expect^{l}_{N}
\big(
h_1(\uzet)
\log h_1(\uzet)
\big|
\zeta_1
\big)
\big)
\\[5pt]
\label{eq:entropydecomp}
&&
+
\expect^{l,1}_N
\left(
\expect^l_N\big(h(\uzet) \big|\zeta_1\big)
\log \expect^l_N\big(h(\uzet) \big|\zeta_1\big)
\right),
\eeq
where in the first term of the right hand side 
\[
h_1(\uzet)
:=
\frac
{h(\uzet)}
{\expect^l_N \big( h(\uzet) \big| \zeta_1 \big)}.
\]

First we bound the first term on the right hand side of
(\ref{eq:entropydecomp}). 
By the induction hypothesis
\beq
\notag
&&
\hskip-2cm
\expect^{l,1}_N
\big(
\expect^l_N\big(h(\uzet) \big| \zeta_1 \big)
\expect^{l}_{N}
\big(
h_1(\uzet)
\log h_1(\uzet)
\big|
\zeta_1
\big)
\big)
\\[5pt]
\notag
\hskip2cm
&=&
\expect^{l,1}_N
\big(
\expect^l_N\big(h(\uzet) \big| \zeta_1 \big)
\expect^{l-1}_{N^{\zeta_1}}
\big(
h_1(\uzet)
\log h_1(\uzet)
\big)
\big)
\\[5pt]
\notag
\hskip2cm
&\le&
W(l-1)
\expect^{l,1}_N
\big(
\expect^l_N\big(h(\uzet) \big| \zeta_1 \big)
D^{l-1}_{N^{\zeta_1}}
\big(
\sqrt {h_1}
\big)
\big)
\\[5pt]
\notag
\hskip2cm
&=&
W(l-1)
\expect^{l,1}_N
\big(
D^{l-1}_{N^{\zeta_1}}
\big(
\sqrt {h(\zeta_1,\cdot)}
\big)
\big)
\\[5pt]
\label{eq:firstterm}
\hskip2cm
&\le&
W(l-1)
D^{l}_{N}
\big(
\sqrt {h}
\big).
\eeq

Next we turn to the second term on the right hand side of
(\ref{eq:entropydecomp}). 
In order to simplify notation in the next argument we denote 
\beq
\label{eq:notation}
\varrho_\alpha:=\frac{N_\alpha}{l},
\qquad
q_\alpha(j):=
\expect^l_N\big(h(\uzet) \ind_{\{\zeta_j=\alpha\}} \big).
\eeq
It is straightforward that for any $K<\infty$  
there exists a finite constant $C=C(K)$ such
that  for any $v\in[0,K]$
\beqs
v\log v \le 
(v-1) +  C \big(\sqrt v -1 \big)^2
\eeqs
and, furthermore, the constant $C$ can be chosen so that 
for any $v>K$
\beqs
v\log v \le 
C v^{3/2}.
\eeqs
Hence, with the notation introduced in (\ref{eq:notation}), 
we get the following upper bound for 
the second term on the right hand side of
(\ref{eq:entropydecomp})
\beq
\label{eq:upper1}
&&
\hskip-6mm
\expect^{l,1}_N
\left(
\expect^l_N\big(h(\uzet) \big|\zeta_1\big)
\log \expect^l_N\big(h(\uzet) \big|\zeta_1\big)
\right)
=
\sum_{\alpha=1}^r
\varrho_\alpha 
\frac{q_\alpha(1)}{\varrho_\alpha} 
\log \frac{q_\alpha(1)}{\varrho_\alpha}
\\[5pt]
\notag
&&
\le 
C
\sum_{\alpha=1}^r
\varrho_\alpha 
\left\{
\left(\sqrt{\frac{q_\alpha(1)}{\varrho_\alpha}} - 1\right)^{\! 2} 
\ind_{\{\frac{q_\alpha(1)}{\varrho_\alpha}\le K\}}
+
\left(\frac{q_\alpha(1)}{\varrho_\alpha}\right)^{\! 3/2} 
\ind_{\{\frac{q_\alpha(1)}{\varrho_\alpha}> K\}}
\right\}.
\eeq
We use the straightforward inequality
\[
\sum_{\alpha=1}^r
\varrho_\alpha 
\left(
\frac{q_\alpha(1)}{\varrho_\alpha}-1
\right)
\ind_{\{\frac{q_\alpha(1)}{\varrho_\alpha}\le K\}}
\le 0. 
\] 
We choose $K$ sufficiently large in order that Lemma 4.1 of
\cite{yau2} can be applied to 
$\{1,2,\dots,l\}\ni j\mapsto \sqrt{q_\alpha(j)/\varrho_\alpha}$. 
Thus we obtain the upper bound
\beq
\notag
\left(\sqrt{\frac{q_\alpha(1)}{\varrho_\alpha}} - 1\right)^{\! 2} 
\ind_{\{\frac{q_\alpha(1)}{\varrho_\alpha}\le K\}}
+
\left(\frac{q_\alpha(1)}{\varrho_\alpha}\right)^{\! 3/2} 
\ind_{\{\frac{q_\alpha(1)}{\varrho_\alpha}> K\}}
&&
\\[5pt]
\label{eq:upper2}
\le 
C' l 
\sum_{j=1}^{l-1} 
\left(
\sqrt{\frac{q_\alpha(j+1)}{\varrho_\alpha}} 
- 
\sqrt{\frac{q_\alpha(j)}{\varrho_\alpha}} 
\right)^{\! 2}. 
&&
\eeq
Putting together (\ref{eq:upper1}) and   (\ref{eq:upper2})
 and returning to the explicit notation we obtain the following upper
 bound for the second term on the right hand side of
 (\ref{eq:entropydecomp}): 
\beq
\notag
&&
\hskip-7mm
\expect^{l,1}_N
\left(
\expect^l_N\big(h(\uzet) \big|\zeta_1\big)
\log \expect^l_N\big(h(\uzet) \big|\zeta_1\big)
\right)
\\[5pt]
\notag
&&
\le 
C'' l 
\sum_{j=1}^{l-1}
\sum_{\alpha=1}^r
\left(
\sqrt{
\expect^l_N\big(h(\uzet)\ind_{\{\zeta_{j+1}=\alpha\}}\big)
}
-
\sqrt{
\expect^l_N\big(h(\uzet)\ind_{\{\zeta_{j}=\alpha\}}\big)
}
\right)^{\! 2}
\\[5pt]
\notag
&&
=
C'' l 
\sum_{j=1}^{l-1}
\sum_{\alpha=1}^r
\left(
\sqrt{
\expect^l_N\big(h(\Theta_{j,j+1}\uzet)\ind_{\{\zeta_{j}=\alpha\}}\big)
}
-
\sqrt{
\expect^l_N\big(h(\uzet)\ind_{\{\zeta_{j}=\alpha\}}\big)
}
\right)^{\! 2}
\\[5pt]
\notag
&&
=
C'' l 
\sum_{j=1}^{l-1}
\left(
\sqrt{
\expect^l_N\big(h(\Theta_{j,j+1}\uzet)\big)
}
-
\sqrt{
\expect^l_N\big(h(\uzet)\big)
}
\right)^{\! 2}
\\[5pt]
\label{eq:secondterm}
&&
\le 
C'' l 
\sum_{j=1}^{l-1}
\expect^l_N
\left(
\left(
\sqrt{h(\Theta_{j,j+1}\uzet)}
-
\sqrt{h(\uzet)}
\right)^2
\right)
=
C'' l 
D^l_N\left(\sqrt h\right).
\eeq
In the second step we used \emph{exchangeability} of the canonical
measures $\pi^l_N$. In the last inequality we note that the
map 
\[
\R_+\times\R_+\ni(x,y)\mapsto\left(\sqrt x-\sqrt y\right)^2
\]
is \emph{convex} and we use  Jensen's inequality.

From (\ref{eq:entropydecomp}), (\ref{eq:firstterm}) and
(\ref{eq:secondterm}) eventually we obtain
\[
W(l)\le W(l-1) + C'' l, 
\]
which yields (\ref{eq:lsigen}).
\end{proof}

\subsection{An elementary probability lemma}
\label{subs:epl}

The contents of the present subsection, in paericular 
Lemma \ref{lemma:epl} and its  Corollary \ref{corollary:epl}  
are borrowed form \cite{tothvalko2}. For their  
proofs see that paper. 

Let
$(\Omega,\pi)$
be a finite probability space and
$\omega_i$, $i\in\Z$
i.i.d.
$\Omega$-valued
random variables with distribution
$\pi$.
Further on let
\beqs
\begin{array}{ll}
\vxi:\Omega\to\R^d,
\quad
&
\vxi_i:=\vxi(\omega_i),
\\[5pt]
\ups:\Omega^{m}\to\R,
\quad
&
\ups_i:=\ups(\omega_{i}\dots,\omega_{i+m-1}).
\end{array}
\eeqs
For
$\vx\in\text{co}(\text{Ran}(\vxi))$
denote
\[
\Ups(\vx):=
\frac
{\expect_\pi\big(\ups_1\exp\{\sum_{i=1}^{m}\vlam\cdot\vxi_i\}\}\big)}
{\expect_\pi\big(\exp\{\vlam\cdot\vxi_1\}\}\big)^{m}},
\]
where $\text{co}(\cdot)$ stands for `convex hull' and 
$\vlam\in\R^d$
is chosen so that
\[
\frac
{\expect_\pi\big(\vxi_1\exp\{\vlam\cdot\vxi_1\}\}\big)}
{\expect_\pi\big(\exp\{\vlam\cdot\vxi_1\}\}\big)}
=
\vx.
\]
For
$l\in\N$
we denote 
\emph{plain} 
block averages by
\[
\overline{\vxi}_l:=\frac1l\sum_{j=1}^l\vxi_j.
\]
Finally, let
$b:[0,1]\to\R$
be a fixed smooth function and denote
\[
M(b):=\int_0^1 b(s)\,ds.
\]
We also define  the block averages 
\emph{weighted by $b$} as
\[
\langle b\,,\,\vxi\rangle_l
:=
\frac1l\sum_{j=0}^l
b(j/l)\vxi_j,
\quad
\langle b\,,\,\ups\rangle_l
:=
\frac1l\sum_{j=0}^l
b(j/l)\ups_j,
\]

The following lemma  relies on elementary probability
arguments:

\begin{lemma}
\label{lemma:epl}
There exists a constant
$C<\infty$,
depending only on $m$, on the joint distribution of
$(\ups_i,\vxi_i)$
and on the function
$b$,
such that the following bounds hold
uniformly in
$l\in\N$
and
$\vx\in({\rm{Ran}}(\vxi)+\dots+{\rm{Ran}}(\vxi))/l$:
\\
(i)
If
$M(b)=0$,
then
\beq
\label{eq:epl11}
\expect
\Big(
\exp
\big\{
\gamma \sqrt{l}
\langle b\,,\,\ups\rangle_l
\big\}
\,\Big|\,
\overline{\vxi}_l=\vx
\Big)
\le
\exp\{C(\gamma^2+\gamma/\sqrt l)\}.
\eeq
(ii)
If
$M(b)=1$
then
\beq
\label{eq:epl12}
\expect
\Big(
\exp
\big\{
\gamma\sqrt{l}
\big(
{\langle b\,,\,\ups\rangle}_l
-
\Ups(\langle b\,,\,\vxi\rangle_l)
\big)
\big\}
\,\Big|\,
\overline{\vxi}_l=\vx
\Big)
\le
\exp\{C(\gamma^2+\gamma/\sqrt l)\}.
\eeq
\end{lemma}

The proof of this lemma appears in \cite{tothvalko2}.

\begin{corollary}
\label{corollary:epl}
There exists a
$\gamma_0>0$,
depending only on $m$, on the joint distribution of
$(\ups_i,\vxi_i)$
and on the function
$b$,
such that the following bounds hold
uniformly in
$l\in\N$
and
$\vx\in({\rm{Ran}}(\vxi)+\dots+{\rm{Ran}}(\vxi))/l$:
\\
(i)
If
$M(b)=0$,
then
\beq
\label{eq:epl1}
\expect
\Big(
\exp
\big\{
\gamma_0 l
\langle b\,,\,\ups\rangle_l^2
\big\}
\,\Big|\,
\overline{\vxi}_l=\vx
\Big)
\le
\sqrt 2.
\eeq
(ii)
If
$M(b)=1$
then
\beq
\label{eq:epl2}
\expect
\Big(
\exp
\big\{
\gamma_0 l
\big(
{\langle b\,,\,\ups\rangle}_l
-
\Ups(\langle b\,,\,\vxi\rangle_l)
\big)^2
\big\}
\,\Big|\,
\overline{\vxi}_l=\vx
\Big)
\le
\sqrt 2.
\eeq
\end{corollary}

\begin{proof}
The bounds
(\ref{eq:epl1})
and 
(\ref{eq:epl2}) 
follow from 
(\ref{eq:epl11}),
respectively, 
(\ref{eq:epl12}) 
by exponential Gaussian averaging.  
\end{proof}

\subsection{Proof of the a priori bounds (Proposition
  \ref{propo:apriori})} 
\label{subs:aprioriproof}

\subsubsection{Proof of the  block replacement bound
  (\ref{eq:obrepl})}

We note first that by simple numerical approximation (no probability
bounds involved)
\beqs
\abs{
\int_{\T}
\abs{
\wih{\ups}(x)
-
\Ups(\wih{\vxi}(x))
}^2
\,dx
-
\frac1n\sum_{j=1}^n
\abs{
\wih{\ups}(j/n)
-
\Ups(\wih{\vxi}(j/n))
}^2
}
\phantom{MMMM}
&&
\\[5pt]
\le
Cl^{-2}
=
o\left(\frac{l^2}{n^2\sigma}\right).
&&
\eeqs
We apply Lemma 
\ref{lemma:varadhan} 
with
\[
{\cal V}_j=
\abs{
\wih{\ups}(j/n)
-
\Ups(\wih{\vxi}(j/n))
}^2.
\]
We use  the bound 
(\ref{eq:epl2}) 
of Lemma 
\ref{corollary:epl}
with the function
$b=a$
of 
(\ref{eq:blocks}).
Note that
$\gamma=\gamma_0l$
can be chosen in 
(\ref{eq:varadhan}). 
This yields the bound
(\ref{eq:obrepl}).

\subsubsection{Proof of the  gradient bound
  (\ref{eq:gradbound})}

Again, we start with numerical approximation:
\beqs
\abs{
\int_{\T}
\abs{
\px\wih{\ups}(x)
}^2
\,dx
-
\frac1n\sum_{j=1}^n
\abs{
\px\wih{\ups}(j/n)
}^2
}
\le
C\frac{n^2}{l^4}
=
o(\sigma^{-1}).
\eeqs
We apply Lemma 
\ref{lemma:varadhan}
 with
\[
{\cal V}_j=
\abs{
\px\wih{\ups}(j/n)
}^2.
\]
We use now the bound 
(\ref{eq:epl1}) 
of Lemma 
\ref{corollary:epl}
with the function
$b=a'$,
where
$a$
is the weighting function from
(\ref{eq:blocks}).
The same choice
$\gamma=\gamma_0l$
applies.  This will yield the bound 
(\ref{eq:gradbound}).


\bigskip

\noindent 
{\bf\large Acknowledgement:} 
We thank the kind
hospitality of Institut Henri Poincar\'e where part of this work
was done. We also acknowledge the financial support of  the
Hungarian Science Foundation (OTKA), grants T26176 and T037685.


\vskip1cm

\hbox{\sc
\vbox{\noindent
\hsize66mm
J\'ozsef Fritz\\
Institute of Mathematics\\
Technical University Budapest\\
Egry J\'ozsef u. 1.\\
H-1111 Budapest, Hungary\\
{\tt jofri{@}math.bme.hu}
}
\vbox{\noindent
\hsize66mm
B\'alint T\'oth\\
Institute of Mathematics\\
Technical University Budapest\\
Egry J\'ozsef u. 1.\\
H-1111 Budapest, Hungary\\
{\tt balint{@}math.bme.hu}
}
}

\end{document}